\documentclass[leqno]{amsart}

\usepackage{amsmath,amsthm,amsfonts,amssymb}
\usepackage[dvips]{epsfig}

\usepackage[colorlinks=true]{hyperref}

\theoremstyle{plain}
\newtheorem{theorem}{Theorem}[section]
\newtheorem{definition}[theorem]{Definition}

\newtheorem{assumption}[theorem]{Assumption}
\newtheorem{lemma}[theorem]{Lemma}
\newtheorem{corollary}[theorem]{Corollary}
\newtheorem{proposition}[theorem]{Proposition}

\theoremstyle{remark}
\newtheorem{remark}[theorem]{Remark}
\newtheorem{example}[theorem]{Example}

\def\la{\left\lvert}
\def\lA{\left\lVert}
\def\ra{\right\rvert}
\def\rA{\right\rVert}

\def\k{{\kappa}}
\def\C{{\mathbb C}}
\def\R{{\mathbb R}}
\def\N{{\mathbb N}}
\def\Z{{\mathbb Z}}
\def\T{{\mathbb T}}
\def\Q{{\mathbb Q}}
\def\A{{\mathbf A}}
 
\def\virgp{\raise 2pt\hbox{,}}

\def\({\left(}
\def\){\right)}
\def\<{\left\langle}
\def\>{\right\rangle}
\def\le{\leqslant}
\def\ge{\geqslant}

\def\d{{\partial}}

\def\eps{\varepsilon}
\def\l{\lambda}
\def\om{\omega}
\def\si{{\sigma}}
\def\F{\mathcal F}
\def\O{\mathcal O}

\DeclareMathOperator{\RE}{Re}
\DeclareMathOperator{\IM}{Im}

\numberwithin{equation}{section}

\begin{document}

%\begin{frontmatter}%
\title[Multiphase geometric optics for NLS]{Multiphase weakly
  nonlinear geometric optics for Schr\"odinger equations}   
\author[R. Carles]{R{\'e}mi Carles}
\address[R. Carles]{Univ. Montpellier~2\\ Math\'ematiques\\
   CC051\\ F-34095~Montpellier}
\address{CNRS, UMR 5149\\ Montpellier\\ F-34095}
\email{Remi.Carles@math.cnrs.fr}
\author[E.~Dumas]{Eric Dumas}
\address[E.~Dumas]{Univ. Grenoble 1\\ Institut Fourier\\ 100, rue des
  Math\'ematiques-BP~74\\ 38402 Saint Martin d'H\`eres cedex\\ France} 
\email{Eric.Dumas@ujf-grenoble.fr}
\author[C. Sparber]{Christof Sparber}
\address[C. Sparber]{Department of Applied Mathematics and Theoretical
  Physics\\ 
  CMS, Wilberforce Road\\ Cambridge CB3 0WA\\ England}
\email{c.sparber@damtp.cam.ac.uk}
\thanks{This work was supported by the French ANR project
  R.A.S. (ANR-08-JCJC-0124-01) and by Award No. KUK-I1-007-43, made by
  King Abdullah University of Science and Technology (KAUST)}  
\begin{abstract}
  We describe and rigorously justify the nonlinear interaction of
  highly oscillatory waves  
  in nonlinear Schr\"odinger equations, posed on Euclidean space or on
  the torus.  
  Our scaling corresponds to a weakly nonlinear regime where the
  nonlinearity affects the leading 
  order amplitude of the solution, but does not alter the 
  rapid oscillations. We consider initial states which are
  superpositions of slowly  
  modulated plane waves, and use the
  framework of Wiener algebras. 
  A detailed analysis of the corresponding nonlinear wave mixing
  phenomena is given, including a geometric interpretation  
  on the resonance structure for cubic nonlinearities.
  As an application, we recover and
  extend some instability results for the nonlinear Schr\"odinger
  equation on the torus in negative order Sobolev spaces. 
\end{abstract}
\maketitle

%\tableofcontents

\section{Introduction}
\label{sec:intro}

\subsection{Physical motivation} The (cubic) \emph{nonlinear
  Schr\"odinger equation} (NLS) 
\begin{equation}
  \label{eq:cubicnls}
  i \d_t u +\frac{1}{2}\Delta u = \lambda 
  |u|^{2}u,
\end{equation}
with $\l\in \R^*$, is one of the most important models in nonlinear science. 
It describes a large number of physical phenomena in 
nonlinear optics, quantum superfluids, plasma physics or 
water waves, see e.g. \cite{SuSu} for a general overview. Independent
of its physical context  
one should think of \eqref{eq:cubicnls} as a description of 
nonlinear waves propagating in a dispersive medium. In the present
work we are interested in 
describing the possible resonant interactions of such waves, often
referred to as \emph{wave mixing}.  
The study of this nonlinear phenomena is of significant mathematical
and physical interest: for example, in the context of fiber optics,
where \eqref{eq:cubicnls}  
describes the time-evolution of the (complex-valued) electric field
amplitude of an optical pulse,  
it is known that the dominant nonlinear process limiting the information 
capacity of each individual channel is given by four-wave mixing,
cf. \cite{HJRB, ZwMe}. Due to its cubic nonlinearity,  
 \eqref{eq:cubicnls} seems to be a natural candidate for the
 investigation of this particular wave mixing phenomena. 
Similarly, four wave mixing appears in the context of plasma physics
where NLS type models are used to describe  
the propagation of Alfv\'en waves \cite{RaTa}. 
Moreover, recent physical experiments have shown the possibility of
matter-wave mixing in \emph{Bose--Einstein condensates} \cite{Deng}.   
A formal theoretical treatment, based on the 
Gross--Pitaevskii equation (\emph{i.e.} a cubic NLS describing
the condensate wave function in a mean-field limit), can be found  
in \cite{TBJ, InTr}. Finally, we also want to mention the closely related
studies on so-called \emph{auto-resonant solutions} of NLS given in
\cite{FrSh} where  
again wave mixing phenomena are used as a method of excitation and
control of multi-phase waves.  

Due to the high complexity of the problem most of the aforementioned
works are restricted to the study of small  
amplitude waves, representing, in some sense, the lowest order
nonlinear effects in  
systems which can approximately be described by a linear superposition
of waves. In addition a slowly varying amplitude approximation is
usually deployed.  
By doing so one restricts himself to resonance phenomena 
which are \emph{adiabatically stable} over large space- and
time-scales. We shall follow this approach by introducing a  
small parameter $0<\eps \ll 1$, which represents the
microscopic/macroscopic scale ratio, and consider  
a rescaled version of \eqref{eq:cubicnls}:
\begin{equation}
  \label{eq:seminls}
  i\eps \d_t u^\eps +\frac{\eps^2}{2}\Delta u^\eps = \lambda \eps
  |u^\eps|^{2}u^\eps.
\end{equation}
This is a \emph{semi-classically scaled} NLS \cite{CaBook}
representing the time evolution of the  
wave field $u^\eps(t,x)$ on macroscopic length- and time-scales. In
the following we 
seek an asymptotic description of  
$u^\eps $ as $\eps\to 0$ on space/time-intervals, which are
independent of $\eps$. 
Note that due to the small parameter $\eps$ in front of the
nonlinearity, we consider 
a \emph{weakly nonlinear regime}. This means that the nonlinearity
does not affect 
the geometry of the propagation, see \S\ref{sec:formel}
below. Technically, it does not show up in 
the eikonal equation, but only in the transport equations determining the
modulation of the leading order amplitudes. In view of these remarks,
the sign of $\l$ (focusing or 
defocusing nonlinearity) turns out to be irrelevant.

\subsection{A general formal computation}
\label{sec:formel}
In order to describe the appearance of the wave mixing in solutions to
\eqref{eq:nls}, we follow the  
\emph{Wentzel-Krammers-Brillouin} (WKB) approach, as first rigorously
settled by Lax \cite{Lax57}. Consider approximate 
solutions of \eqref{eq:seminls} in the form of high-frequency wave
packets, such as  
\begin{equation}
  \label{eq:wkb}
  a(t,x)e^{i \phi(t,x)/\eps}.
\end{equation}
For such a single mode to be an approximate solution, it is necessary
that the rapid oscillations are  
carried by a phase $\phi$ which solves the eikonal equation (see
\cite{CaBook}, where  
also other regimes, in terms of the size of the coupling constant, are
discussed):
\begin{equation}
  \label{eq:eikonal}
  \d_t \phi+\frac{1}{2}\lvert \nabla \phi\rvert^2 =0.
\end{equation}
Nonlinear interactions of high frequency waves are then found by
considering superpositions of wave packets \eqref{eq:wkb}.  
By the cubic interaction, three phases $\phi_1$, $\phi_2$ and $\phi_3$
generate 
\begin{equation*}
 \label{eq:phicreation}
 \phi= \phi_1-\phi_2+\phi_3.
\end{equation*}   
The corresponding term is 
relevant at leading order if and only if this new phase $\phi$ is
\emph{characteristic}, \emph{i.e.} solves the 
eikonal equation \eqref{eq:eikonal} while also each $\phi_j$, $j
=1,2,3$ does so.  
More generally, we will have to construct a set of phases $\{ \phi_j \}_{j
  \in J}$, for some index set $J \subset \Z$, such that each $\phi_j$
is characteristic,  
and the set is \emph{stable} under the nonlinear interaction. That is,
if $k, \ell, m \in J$ are such that $\phi = \phi_k-\phi_\ell+\phi_m$ is
characteristic, then $\phi \in \{ \phi_j \}_{j \in J}$.  
Given some index $j \in J$, the set of (four-wave) resonances leading
to the phase $\phi_j$ is then   
\begin{equation*}
 I_j = \big \{  (k,\ell,m) \in J^3 \ ;\   \phi_k-\phi_\ell+\phi_m = \phi_j  \big \} .
 \label{eq:resonance}
\end{equation*}
%\begin{remark}
%In view of the usual microlocal theory, we tend to think that a phase
%should be a smooth function $\phi$, solution to \eqref{eq:eikonal},
%with a  
%gradient $\nabla\phi$ that does not vanish on any open set. However,
%since nonlinear interactions may create non-oscillating modes, we will
%also consider $\phi\equiv0$ to be an admissible phase. 
%\end{remark}
One of the tasks of this work to 
study the structure of $I_j$. 
A first important step is obtained by plugging $\phi
= \phi_k-\phi_\ell+\phi_m $ into \eqref{eq:eikonal}, since then, an
easy 
calculation shows that $\phi$ is characteristic  if and only if the following 
\emph{resonance condition} is satisfied:
\begin{equation}
  \label{eq:resonance3}
  \(\nabla \phi_\ell -\nabla\phi_m\)\cdot \(\nabla \phi_\ell -\nabla \phi_k\)=0. 
 \end{equation}
Obviously this is a quite severe restriction in one spatial dimension,
while in higher dimensions there are many possibilities to satisfy
\eqref{eq:resonance3}. 
In order to gain more insight we shall restrict ourselves from now on
to the case of \emph{plane waves} (\emph{i.e.} linear phases,
see \S\ref{sec:phasegeneral}),  
This choice allows for a more detailed mathematical study and is also
the most important case from the physical point of view,
cf. \cite{TBJ, FrSh}. 
The precise mathematical setting is then as follows.

\subsection{Basic mathematical setting and outline}
 
In the following the space variable $x\in \mathcal M$ will either
belong to the whole Euclidean space $\mathcal M=\R^d$, or to the torus
$\mathcal M=\T^d$ (we denote 
$\T=\R/2\pi\Z$), for some $d\in \N$. The latter can be motivated by
the fact that numerical simulations of \eqref{eq:nls}  
are mainly based on pseudo-spectral schemes and thus naturally posed
on $\T^d$, see e.g. \cite{BSM1, BSM2}.  We then consider the initial
value problem for the  slightly more general NLS   
\begin{equation} 
  \label{eq:nls}
  i\eps \d_t u^\eps +\frac{\eps^2}{2}\Delta u^\eps =  \, \lambda \eps
  |u^\eps|^{2\si}u^\eps \quad ; \quad 
  u^\eps (0,x) =   u_0^\eps(x),
\end{equation}
where $\si\in \N^*$. Although we obtain the most precise results 
(concerning the geometry of resonances, in particular) in the case of 
the cubic nonlinearity ($\si=1$), we are in fact able to rigorously 
justify WKB asymptotics also for higher order nonlinearities.  
We assume that \eqref{eq:nls} is subject to an initial data $u_0^\eps$, which is 
assumed to be close (in a sense to be made precise in \S\ref{sec:justif}) to superposition of highly
oscillatory plane waves, \emph{i.e.}
\begin{equation}\label{eq:CIgen}
  u_0^\eps(x) \approx \sum_{j\in
  J_0}\alpha_j\(x\)e^{i \kappa_j \cdot x /\eps}, 
\end{equation}
where $J_0\subseteq \Z$ is a (not necessarily finite) given index set. 
In the Euclidean case we allow for wave vectors $\kappa_j \in \R^d$, 
whereas on $\mathcal M= \mathbb T^d$ we impose $\kappa_j \in \Z^d$.
Moreover, in the latter case, we choose $\alpha_j$ to be independent
of $x \in \T^d$, so that
\eqref{eq:CIgen} corresponds to an expansion in terms of Fourier series
(with $\eps^{-1}\in\N$). The case of $x$-dependent $\alpha_j$'s on
$\T^d$ could be considered as well, by reproducing the analysis on
$\R^d$. We choose not to do so here, since it brings no real new 
information. 

In particular, for $x\in \T$ (the one-dimensional torus), our analysis 
leads to a remarkably simple approximation. 
\begin{theorem}\label{theo:torus}
  For $x\in \T$, consider \eqref{eq:nls} with $\si=1$. Suppose that
  the initial data are of the form  \eqref{eq:CIgen} with $\k_j \equiv
  j \in \Z$ and $(\alpha_j)_j\in \ell^1(\Z)$. 
  
  Then for all $T>0$,
  there exist $C=C(T)$ and $\eps_0>0$, such that for all $\eps \in
  ]0,\eps_0]$, with $1/\eps\in\N^*$, it holds 
\begin{equation*}
  \sup_{t\in [0,T]}\left\lVert
    u^\eps(t)-u_{\rm app}^\eps(t)\right\rVert_{L^\infty(\T)}\le C\eps,
\end{equation*}
where the approximate solution $u_{\rm app}^\eps$ is given by
\begin{equation*}
  u_{\rm app}^\eps(t,x) = \sum_{j\in\Z}
  \alpha_j e^{-i\l t(2M - |\alpha_j|^2)} e^{i(jx-\frac{1}{2}j^2t)/\eps},
\text{ and } M = \sum_{k\in
    \Z}|\alpha_k|^2.  
\end{equation*}
\end{theorem}
We see that at leading order, the nonlinear interaction shows up
through an explicit modulation at scale $\O(1)$. It is well known that
the one-dimensional 
cubic Schr\"odinger equation is completely integrable
(see \cite{Its76,MA81} for the periodic case). However, this aspect
does not play any role in the proof of Theorem~\ref{theo:torus}, which in itself
can be seen as a consequence of the more general result stated in
Theorem~\ref{theo:justif}. 
On the other hand, several aspects in the discussion on possible phase
resonances and the creation of  amplitudes seem to be
specific to both properties $d=1$ and $\si=1$ (see \S\ref{sec:phase}
and \S\ref{sec:profile}). 

In order to prove Theorem~\ref{theo:justif}, and henceforth also
Theorem \ref{theo:torus}, we need to set up a rigorous  
multiphase WKB approximation for solutions to  \eqref{eq:nls}. To this
end, there are essentially two steps  
needed in our analysis. First, we detail the
approach sketched above by examining the possible resonances between the 
phases, and analyzing the evolution and/or the creation of
the corresponding profiles $a_j$. The second step then consists in making
this approach rigorous: we construct the profiles $a_j$, and 
show that the obtained \emph{ansatz} is a satisfactory approximation
of the exact solution $u^\eps$, up to $\O(\eps)$ in a space contained
in $L^\infty(\mathcal M)$. As it is standard, we prove in fact a
stronger stability result: Starting from any approximate solution
$u^\eps_{\rm app}$ constructed on profiles, we show that, for any initial
data close (as $\eps$ goes to zero) to $u^\eps_{\rm app\mid t=0}$,
there exists an exact solution which is close to $u^\eps_{\rm app}$,
on some time interval independent of $\eps$ (which, for $\eps$ small
enough, may be chosen as any finite time up to which $u^\eps_{\rm
  app}$ is defined).  

In the case of a single oscillation only, it suffices to multiply
$u^\eps$ by $e^{-i\phi/\eps}$ to filter out rapid oscillations, see
\cite{CaBook}. In the case where 
several phases are present, this strategy obviously fails. 
To overcome this issue, a fairly general mathematical approach, which has
proved efficient in several contexts (see
e.g. \cite{Gues93,RauchUtah,RauchNotes}), consists in working in rescaled 
Sobolev spaces, usually denoted by $H^s_\eps$, for $s>0$. These are
the usual Sobolev spaces, where derivatives are scaled by
$\eps$, in order to account for the spatial oscillations at scale
$\eps$. More precisely, if  $s\in \N$,
\begin{equation*}
  \|f\|^2_{H^s_\eps}:=\sum_{|\alpha|\le s} \|(\eps\d)^\alpha
    f\|_{L^2}^2. 
\end{equation*}
However, due to the negative power of $\eps$ in the associated
Gagliardo--Nirenberg inequalities, this technique usually demands to
construct approximate solution with a high order of precision (see
\cite{GMS} for a closely related study on the interaction of
high-frequency waves in periodic potentials).  
% We sketch the application of this technique in an appendix.  
Another, more sophisticated, approach consists in filtering out the
rapid oscillations in terms of the free evolution group, as in
\cite{Schochet}. In the present work though, we shall use a simpler
approach, which allows us to justify the multiphase weakly nonlinear
WKB analysis in a remarkably straightforward way. This approach relies
on the use of Wiener algebras, as introduced in \cite{JMRWiener}, and
further developed in \cite{Lan01,BL02,MColinLannes}. This analytical
framework is particularly convenient in the case of plane waves, but
could probably be extended to more general situations, up to some geometric
constraints on the phases. However, the first step of the analysis,
\emph{i.e.} describing all possible resonances,  
becomes much more intricate, see e.g. \cite{JMRENS,JMRDuke}. 

As well shall see during the course of the proof, the use of Wiener
algebras has several advantages 
on the technical level. We point out that this framework makes it possible
to justify the WKB approximation with an error estimate of order
$\O(\eps)$ without constructing correctors (which would have to be of
order $\eps$ or even smaller, when working in $H^s_\eps$ spaces, see
e.g. \cite{Gues93,RauchUtah}, or \cite{CaMaSp04,GMS} in the NLS case). 

\subsection{An application to instability}

As an application of the semi-classical analysis for \eqref{eq:nls},
we recover the main result in \cite{CCTper} (see also \cite{BGTMRL}),
concerning NLS in the periodic case. This result has been established in
the case $d=1$, and is hereby extend to higher dimensions. We
also propose a variation on a result in \cite{Molinet} (see assertion
3 in the theorem below).  
\begin{theorem}\label{theo:instab}
  Let $d\ge 1$, $\si \in \N^*$ and $\lambda \in \{\pm 1\}$. Fix
  $s<0$. 
  \smallbreak

\noindent $1.$ For all $\rho>0$, we 
  can find a solution $u$ to 
  \begin{equation}\label{eq:nlsT}
    i\d_t u + \frac{1}{2}\Delta u= \lambda
    |u|^{2\si}u,\quad x\in \T^d, 
  \end{equation}
with $\|u(0)\|_{H^s(\T^d)}<\rho$, such that for all $\delta>0$,
there exists $\widetilde u$ solution to \eqref{eq:nlsT} with
\begin{equation*}
  \lVert u(0) - \widetilde u(0)\rVert_{H^s(\T^d)}<\delta,
\end{equation*}
and
\begin{equation*}
 \sup_{0\le t\le \delta} \left\lvert \int_{\T^d} \(u(t,x)-\widetilde
    u(t,x)\)dx\right\rvert\ge c\rho, 
\end{equation*}
for some constant $c>0$ independent of $\rho$ and $\delta$. In
particular, the solution map fails to be continuous as a map from $H^s(\T^d)$
to $H^k(\T^d)$,  
no matter how close to $-\infty$ the exponent $k$ may be.
\smallbreak
\noindent $2.$ Suppose $\si\ge 2$. For any $\rho > 0$ and $\delta > 0$ 
there exist smooth solutions $u$, $\widetilde u$ of \eqref{eq:nlsT} such that
$u(0)-\widetilde u(0)$ is equal to a constant  
of magnitude at most $\delta$, and 
\begin{equation*}
  \|u(0)\|_{H^s(\T^d)}+\|\widetilde u(0)\|_{H^s(\T^d)}\le \rho\quad
  ;\quad \sup_{0\le 
        t\le \delta} \left\lvert \int_{\T^d} \(u(t,x)-
    \widetilde u(t,x)\)dx\right\rvert\ge c\rho, 
\end{equation*}
for some constant $c>0$ independent of $\rho$ and $\delta$. 
\smallbreak
\noindent $3.$ For any $t\not =0$, the flow-map associated with \eqref{eq:nlsT}
is discontinuous as 
a map from  $L^2(\T^d)$, equipped with its weak topology, into the
space of distributions  $\(C^\infty(\T^d)\)^*$ at any \emph{constant}
$\alpha_0\in \C\setminus\{0\} \subset L^2(\T^d)$. 
\end{theorem}
 We show in
\S\ref{sec:instab} that the above instability result 
can be viewed as a consequence of multiphase weakly nonlinear
geometric optics. 
The first two assertions are an extension of the results in
\cite{CCTper}, so we shall not comment on their meaning, and refer to the
discussion in \cite{CCTper}. We invite the reader to consult
\cite{Molinet} for a stronger instability result in the
one-dimensional case: indeed, when $d=\si=1$, the author shows the third
point in the above statement for any $\alpha_0\in
L^2(\T)\setminus\{0\}$, not necessarily constant.

\subsection{Structure of the paper}
\label{sec:structure}

We first study in detail the case of the cubic nonlinearity ($\si=1$).
In
\S\ref{sec:phase}, we consider the set of resonant phases, and in
\S\ref{sec:profile}, we analyze the corresponding amplitudes. The case
of higher order nonlinearities is treated in \S\ref{sec:higher}. In
\S\ref{sec:analytical}, we set up the analytical framework, with which
a general stability result (of which Theorem~\ref{theo:torus} is a
straightforward consequence) is established in
\S\ref{sec:justif}. Theorem~\ref{theo:instab} is proved in
\S\ref{sec:instab}. Finally, in an appendix, we sketch how the
previous semi-classical analysis can be adapted 
to more general sets of initial plane waves (including generic finite sets
of wave vectors).

\subsubsection*{Acknowledgments} The first author wishes to thank
Thierry Colin and David Lannes for preliminary discussions on this
subject.

\section{Analysis of possible resonances in the cubic case}
\label{sec:phase}

In this section, we show that when $\si=1$, the set of relevant phases
can be described in a fairly detailed way. 

\subsection{General considerations} \label{sec:phasegeneral}
We seek an approximation of the form 
\begin{equation*}\label{eq:approx}
  u^\eps(t,x)\thickapprox \sum_{j\in
  J}a_{j}(t,x)e^{i \phi_j(t,x)/\eps},
\end{equation*}
where here and in the following $J\subset \Z$ denotes the index set of
relevant phases $\phi_j$ (yet to be determined).  
Note that using $J$ is only a renumbering, so that $j\neq k
\Rightarrow \phi_j\neq\phi_k$. In the case  $x\in \T^d$, one simply
drops the dependence of $a_j$ upon 
$x$. In general $J_0 \subsetneq J$, \emph{i.e.} we usually need
to take into account more phases in \eqref{eq:approx} than we are
given initially.  

As a first step we need to determine the characteristic phases
$\phi_j(t,x) \in \R$.  
For plane-wave initial data of the form \eqref{eq:CIgen} we are led to the
following initial value problem  
\begin{equation*}
  \label{eq:eikivp}
   \d_t \phi_j+\frac{1}{2}\lvert \nabla \phi_j \rvert^2 =   0 \quad ;
   \quad \phi_j(0,x) =    \kappa_j \cdot x, 
\end{equation*}
the solution of which is explicitly given by
\begin{equation}\label{eq:planewave}
  \phi_j(t,x) =\kappa_j\cdot x -\frac{t}{2}|\kappa_j|^2.
\end{equation}
Recall that for $x \in \R^d$, we assume $\k_j \in \R^d$, whereas in
the case $x\in \T^d$, we restrict ourselves to $\k_j \in \Z^d$. 
Of course, these phases $\phi_j$ remain smooth for all time,
\emph{i.e.} \emph{no caustic} appears.  

In the cubic case $\si = 1$, the set of resonances leading to the phase 
$\phi_j$ is therefore given by 
\begin{equation*}
\label{eq:planeresonance}
 I_j=\{  (k,\ell,m) \in J^3\ ;\ \kappa_k-\kappa_\ell +\kappa_m
 =\kappa_j,\ |\kappa_k|^2 -|\kappa_\ell|^2 + |\kappa_m|^2 = |\kappa_j|^2\},
\end{equation*}
and the corresponding resonance condition \eqref{eq:resonance3} becomes
\begin{equation}\label{eq:resonanceplane}
  \(\k_\ell -\k_m \)\cdot \(\k_\ell -\k_k\)=0. 
\end{equation}
As we shall see, this condition provides several insights on the
structure of four-wave resonances. 

\subsection{The one-dimensional case}
\label{sec:phases1d}

For $d=1$ the condition 
\eqref{eq:resonanceplane} implies that if $(k,\ell,m)\in I_j$, then 
$\k_\ell= \k_m$, or $\k_\ell= \k_k$. 
Therefore, when $d=1$, the set $I_j$ is fully described by: 
\begin{equation*}\label{eq:1Dresonance}
  I_j =\{ (j, \ell,\ell),\ (\ell,\ell,j)\ ;\ \ell \in J \}, 
\end{equation*}
\emph{no new phase} can be generated by a cubic interaction. 
\smallbreak
In higher dimensions, however, the situation is much more complicated
and heavily depends on the number of initial modes. 

\subsection{Multi-dimensional case  $d\ge2$}

We start with the simplest multiphase situation and proceed from there
to more complicated cases. Eventually we shall arrive at a geometric  
interpretation for the generic case.

\subsubsection{One or two initial modes}
\label{sec:2phases}

If we start from only two initial modes, $ \sharp J_0 = 2$, the
resonance condition \eqref{eq:resonanceplane} implies that the cubic
interaction between  
these two phases \emph{cannot create} a new characteristic phase. In
other words, $u^\eps$ exhibits at most two rapid oscillations at
leading order. Recalling that $\phi=0$ is an
admissible phase, the case of a single
initial phase $ \sharp J_0 = 1$, is therefore included (if one of the initial amplitudes is set equal to zero).  
We want to emphasize that the case of at most two initial phases is rather particular, since
\eqref{eq:resonanceplane} implies that the situation is the same for
all spatial dimensions $d \ge 1$.  

\begin{remark}
In addition, the fact that two phases cannot create a new one 
extends also to higher order (gauge invariant) nonlinearities 
$f(z)=\l|z|^{2\si}z$, for $\si \in \N$, $\si \ge 2$, see \S\ref{sec:higher}.
\end{remark}

\subsubsection{Three or four initial modes}
\label{sec:3phases}

This case can be fully understood by the following geometric insight, 
already noticed in \cite{Iturbulent}:
\begin{lemma} \label{lem:rectangles}
Let $d\ge2$, and $k,\ell,m$ belong to $J$. Then, $(\k_k,\k_\ell,\k_m) \in I_j$ 
precisely when the endpoints of the vectors $\kappa_k, \kappa_\ell, \kappa_m,
\kappa_j$ form \emph{four corners of a non-degenerate rectangle}  
with $\kappa_\ell$ and $\kappa_j$ opposing each other, or when this
quadruplet corresponds to one of the two following degenerate cases:
$(\k_k=\k_j, \k_m=\k_\ell)$,  
or $(\k_k=\k_\ell, \k_m=\k_j)$. 
\end{lemma}
\begin{remark}
In the degenerate cases, no new phase is created. 
\end{remark}
\begin{proof} We recall the argument given in \cite{Iturbulent}, by
  first noting that  
the relations between $(\kappa_j,\kappa_k, \kappa_\ell, \kappa_m)$
formulated in \eqref{eq:planeresonance}, are equivalently fulfilled by  
$(\kappa_k-\kappa, \kappa_\ell-\kappa, \kappa_m-\k,
\kappa_j-\k)$, for any $\k \in \R^d$ (resp. $\kappa\in \mathbb
Z^d$). This is easily seen by expanding the second relation in \eqref{eq:planeresonance}  
and inserting the first one.  
Thus, choosing $\kappa = \kappa_j$,  
it therefore suffices to prove this geometric interpretation for
$\kappa_j = 0$,  which consequently shows: $\k_k + \k_m = \k_\ell$ such that $\kappa_k \cdot \kappa_m = 0$, 
by the law of cosines.
\end{proof}
In summary, we conclude that three initial (plane-wave) phases
\emph{create at most one new phase}, such that the corresponding four
wave vectors form a rectangle. When the initial wave vectors $\{ \k_j
\}_{j \in J_0}$ are chosen such that their endpoints form the
four corners of a rectangle, no new phase can be created by the cubic
nonlinearity and $u^\eps$ exhibits only four rapid oscillations. We
close this subsection with two illustrative examples. 
\begin{example}\label{ex:phase0}
   Let $d=2$. Consider
  $  \k_1=(0,1)$, $\k_2=(1,1)$ and  $\k_3=(1,0)$.
  The cubic interaction creates the zero mode $\phi_4 \equiv 0$.
\end{example}
\begin{example}\label{ex:phasenon0}
  Again let $d=2$, with now $\k_1=(1,1)$, $\k_2=(1,2)$ and
  $\k_3=(3,2)$. 
In this case, we create a non-zero phase $\phi_4$, with corresponding wave
vector $\k_4=(3,1)$. %(see Figure~\ref{fig:ex25}). 
%\begin{figure}
  %\begin{center}
  %\scalebox{0.6}{\input{ex25.pdf}}
%\caption{Example~\ref{ex:phasenon0}}
 % \end{center}
%\label{fig:ex25}
%\end{figure}
\end{example}

The geometric insight gained above then directly leads us to the
following description of the resonant set $I_j$ in the
general case. 

\subsubsection{The general case}
\label{sec:phasegen}

We are given a countable (possibly finite) number of initial phases
$\{ \phi_j \}_{j \in J_0}$  with corresponding wave vectors $\{ \k_j
\}_{j \in J_0}$. From the discussion of the previous paragraph it is
clear that there are two possible situations: 
\begin{itemize}
\item[{(a)}] Either, it is \emph{impossible} to create  a \emph{new}
  rectangle from any possible subset   
$\tilde J_0 \subset J_0$, such that $ \sharp \tilde J_0 = 3$. If so,
then no new phase can be created.  This is the generic case. 
\item[{(b)}] Or, starting from an initial (finite or countable)
  set $S_0=\{\kappa_j\}_{j\in J_0}$, we may obtain a first generation
  $S_1=\{\kappa_j\}_{j\in J_1}$  
with $J_0 \subset J_1$ (\emph{i.e.} $S_0 \subset S_1$) in the
following way: we add to $S_0$ all points $\kappa\in\R^d$, such that
there exist $\tilde{J}_0 \subset J_0$ with $\sharp \tilde{J}_0 = 3$,
and  
such that $\{\k_j\}_{\tilde{J}_0} \cup \{\kappa\}$ is a
rectangle. Note that, if $J_0 \subset \Z^d$, then $J_1 \subset \Z^d$.  
By a recursive scheme, we are led to a (finite or countable) set $S$
which is stable under the completion of right-angled triangles   
formed of points from this set, into
rectangles. Furthermore, if $S_0 \subset \Z^d$, then $S
\subset \Z^d$. 
\end{itemize}

\begin{example} As already seen, the simplest examples for
  possibility (a) are the cases $\sharp  J_0 \le 3$ when the triangle 
formed by the endpoints of the  
  considered wave vectors has no right
angle, or $\sharp  J_0 = 4$, where the four initial phases are chosen
such that their corresponding wave vectors $\{ \k_j \}_{j \in J_0}$
already form the corners of a rectangle.  
\end{example}

From a finite number of initial phases, possibility (b) may lead to a
finite as well as to an infinite set $J$. Even for $d=2$, we have: 

\begin{example}
In the plane $\R^2$, start with 
$$J_0 = \{ (-1,1), (0,1), (0,0), (1,0) \}.$$
The first generation is then 
$$J_1 = \{ (-1,1), (0,1), (1,1), (-1,0), (0,0), (1,0) \} = J_0 \cup \{
(1,1), (-1,0) \},$$ 
and the second one is
$$J_2 = J_1 \cup \{ (0,2), (0,-1) \}.$$ 
One easily sees that this generates $J=\Z^2$.
\end{example}

As a conclusion, the set of phases $\{\phi_j\}_{j\in J}$ may be finite
or infinite, but has the following property. 
\begin{proposition} \label{prop:J}
Let $\si=1$, and consider any triplet of wave vectors from 
$S=\{\kappa_j\}_{j\in J}$. Then, either the corresponding triangle 
has no right angle, or the fourth corner of the associated rectangle 
belongs to $S$.  
\end{proposition}

\section{Analysis of the amplitude system in the cubic case}
\label{sec:profile}

From the previous section, in general we have to expect the generation
of new phases by the four-wave resonance. However, it may happen that
not all of them are actually  present in our approximation \eqref{eq:approx},   
since  the corresponding profile $a_j(t,x)$ has to be non-trivial.  
\smallskip

Indeed, if we plug the ansatz \eqref{eq:approx} into \eqref{eq:nls}
the terms of order $\O(1)$ are identically zero since all the
$\phi_j$'s are characteristic.  
For the $\O(\eps)$ term, we project on the oscillations associated to
$\phi_j$, which yields the following system of transport equations: 
\begin{equation}\label{eq:transportsystem}
  \forall j \in J, \quad \d_t a_j +\k_j \cdot \nabla a_j 
  = -i \l\sum_{(k,\ell,m)\in I_j}a_k \overline
  a_\ell a_m\quad ; \quad a_j(0,x) =    \alpha_j ( x),
\end{equation}
with obviously $\nabla a_ j = 0$ in the case where $x \in \T^d$. In
the following we will perform a qualitative analysis of the system
\eqref{eq:transportsystem}, postponing the rigorous existence and
uniqueness analysis to \S\ref{sec:constrprof}.  Having in mind the
discussion from \S\ref{sec:phase} we distinguish the case $d=1$ from
the case $d\ge 2$.  

\subsection{The case $d=1$}
\label{sec:1d}

Let $j\in J$, and recall that $I_j$ is particularly simple in $d=1$:
\begin{equation*}
  I_j =\{ (j, \ell,\ell),\ (\ell,\ell,j)\ ;\ \ell \in J \}. 
\end{equation*}
Using this, \eqref{eq:transportsystem} simplifies to
\begin{equation}\label{eq:transport1d}
  \(\d_t+\k_j \partial_x\) a_j = -2 i\l \sum_{\ell \in J}|a_\ell|^2
  a_j+i\l|a_j|^2 a_j \quad ; \quad a_j (0,x) =    \alpha_j ( x).
\end{equation}
In particular, the evolution of a zero profile $\alpha_j\equiv 0$ is
necessarily trivial, that is $a_j(t,x)\equiv 0$.
This non-generation of profiles leads to the same conclusion as 
\S\ref{sec:phases1d}: No new mode can be created, if it is not present 
initially (and the reason is the same as in \S\ref{sec:phases1d}: $a_j$ 
factors out in \eqref{eq:transport1d} just because for any $(\ell_1,
\ell_2,\ell_3)\in I_j$, we have $\ell_1=j$ or $\ell_3=j$). We shall see that the
multi-dimensional situation is quite different but 
first examine the situation for $x\in \T$ and $x\in \R$ in
more detail. 

\subsubsection{The case $x\in \T$}
In this case, we readily obtain that $ |a_j|^2$ does not depend on
  time. This is due to the fact that \eqref{eq:transport1d} yields:
  $i\d_t a_j \in \R a_j$ and hence $\d_t |a_j|^2=0$,
  for all $j\in \Z$. In particular we get that 
  \begin{equation*}
    M=\|u^\eps(0)\|_{L^2}^2 =\sum_{j\in J } |\alpha_j|^2
    =\|u^\eps(t)\|_{L^2}^2, \quad \forall t\in \R.
  \end{equation*}
  The conserved quantity $M$ corresponds to the total mass of the
  exact solution $u^\eps$. Using this, we rewrite \eqref{eq:transport1d} as
 \begin{equation*}
  \frac{d}{dt} a_j = -i\l \(2 M-|\alpha_j|^2\) a_j,
\end{equation*}
which yields an explicit formula for the (global in time) solution
\begin{equation*}\label{eq:1dtorusprofile}
  a_j(t)=\alpha_j  e^{-i\l t \(2 M-|\alpha_j |^2\)}. 
\end{equation*}
We observe that in the case of the one-dimensional torus, the
interaction of the profiles $a_j$ is particularly simple. Nonlinear
effects lead to phase-modulations only. 

\subsubsection{The case $x\in \R$}

Here, in contrast to the situation on $\T$, the modulus of $a_j$ is no
longer conserved, since we can only conclude from
\eqref{eq:transport1d} that  
\begin{equation*}
  \(\d_t+\k_j\d_x\) |a_j|^2 =0,
\end{equation*}
and thus
\begin{equation*}
  |a_j(t,x)|^2 =|\alpha_j(x-t\k_j)|^2.
\end{equation*}
In particular we readily see that for all $j\in J$ we have
\begin{equation}\label{eq:L2con}
\|a_j(t)\|_{L^2}= \|\alpha_j\|_{L^2},  \quad \forall t\in \R .
\end{equation}
Moreover, we still have an explicit representation for the solution of
\eqref{eq:transport1d} in the form 
\begin{equation} \label{eq:selfmodulation}
  a_j(t,x)=\alpha_j(x-t\k_j)e^{i S_j(t,x)},
\end{equation}
for some real-valued phase $S_j$, yet to be computed. 
In view of the identity
\begin{equation*}
\(\d_t+\k_j\d_x\)a_j(t,x) = i \alpha_j(x-t\k_j) e^{i\l S_j(t,x)}
\(\d_t+\k_j\d_x\)S_j(t,x),
\end{equation*}
equation~\eqref{eq:transport1d} implies
\begin{align*}
\( \(\d_t+\k_j\d_x\)S_j(t,x) \) \alpha_j(x-t\k_j) = &  
\l \( -2 \sum_{\ell \in J}|\alpha_\ell(x-t\k_\ell)|^2 +
|\alpha_j(x-t\k_j)|^2 \)  \\ 
& \, \times 
\alpha_j(x-t\k_j).
\end{align*}
One easily sees that it is sufficient 
to impose
\begin{equation*}
\d_t\( S_j(t,x+t\k_j) \) = -2 \l \sum_{\ell \in
  J}|\alpha_\ell(x+t(\k_j-\k_\ell))|^2  
+ \l |\alpha_j(x)|^2 ,
\end{equation*}
which yields 
\begin{equation} \label{eq:selfmodulationformula}
  \begin{aligned}
  S_j(t,x) &= 
-2 \l \int_0^t \( \sum_{\ell \in J
  \setminus\{j\}}|\alpha_\ell(x+(\tau-t)\k_j-\tau\k_\ell))|^2  
d\tau \) \\
&\quad - t \l|\alpha_j(x-t\k_j)|^2.
    \end{aligned}
\end{equation}
This formula, together with \eqref{eq:selfmodulation} describes the
modulation of the profile $a_j(t,x)$. As in the case of the torus, amplitudes
are transported linearly. Only the (slow) phases $S_j$ undergo
nonlinear effects, which are more complicated as before but  
still explicitly described in terms of the initial data.

\subsection{The case of one or two modes for $d\ge 1$}
\label{sec:2profiles}

We have already seen in \S\ref{sec:2phases} that the case of two
initial modes is special, since we get a closed system for all
$d\ge1$. Indeed if we start from two phases and two associated
profiles, say $a_j$ and $a_\ell$, the
system~\eqref{eq:transportsystem} simplifies to: 
\begin{equation*}
  \begin{aligned}
    &\d_t a_j + \k_j\cdot\nabla a_j = -i\l \(|a_j|^2+2|a_\ell|^2\)a_j,
    \quad ; \quad a_j(0,x) =    \alpha_j ( x),\\ 
&\d_t a_\ell +\k_\ell\cdot\nabla a_\ell = -i
\l\(2|a_j|^2+|a_\ell|^2\)a_\ell,  \quad \, ; \quad a_\ell (0,x) =
\alpha_\ell ( x). 
  \end{aligned}
\end{equation*}
Note that if initially one of the two profiles is identically zero, it
remains zero for all times and hence, we are back in the  
situation of a usual single-phase WKB approximation. In particular we
compute explicitly for: 
\begin{itemize}
\item Two modes, on $\T^d$:
  \begin{align*}
   a_j(t) = \alpha_j e^{-i \l t (2 |\alpha_\ell |^2 + |\alpha_j
     |^2)}\quad ;\quad 
a_\ell(t) = \alpha_\ell  e^{-i \l t (2 |\alpha_j |^2 + |\alpha_\ell |^2)}.
  \end{align*}
\item Two modes, on $\R^d$:
  \begin{align*}
    a_j(t,x) &= \alpha_j(x-t\k_j) e^{-i\l \(2 \int_0^t
      |\alpha_\ell(x+(\tau-t)\k_j-\tau\k_\ell)|^2  
d\tau + t |\alpha_j(x-t\k_j)|^2\)},\\
a_\ell(t,x) &= \alpha_\ell(x-t\k_\ell) e^{-i\l \(2 \int_0^t
  |\alpha_j(x+(\tau-t)\k_\ell-\tau\k_j)|^2 d\tau + t
  |\alpha_\ell(x-t\k_\ell)|^2\)}.  
  \end{align*}
\end{itemize}
Again, these solutions exhibit (nonlinear) self-modulation of phases only, and 
exist for all times $t \in \R$, a property which is a-priori not clear in the
general case.

\subsection{Creation of new modes when $d\ge 2$}
\label{sec:creationphase}

A basic difference between the one-dimensional case and the
multidimensional situation is that the conservation law
\eqref{eq:L2con} does not remain valid when $d\ge 2$.  
However, we are still able to prove that the total mass is conserved.
\begin{lemma} For any solution to \eqref{eq:transportsystem} it holds
  \begin{equation}\label{eq:masscon}
    \frac{d}{dt}\sum_{j\in J} \|a_j(t)\|_{L^2}^2=0. 
  \end{equation}
  \end{lemma}
\begin{proof} The assertion follows from the more general identity
\begin{equation*}
 \sum_{j\in J} \(\d_t+ \k_j\cdot \nabla\) |a_j|^2 =0,
\end{equation*}
since, by definition we have
\begin{equation*}
  \sum_{j\in J} \(\d_t +\k_j\cdot \nabla\)|a_j|^2 = \IM\( \l 
  \sum_{j\in J} \sum_{(k,\ell,m)\in I_j}\overline a_ja_k\overline
  a_\ell a_m\).
\end{equation*} 
This sum is zero by symmetry, since for each quadruplet $(j,k,\ell,m)\in J^4$ of
indices the quadruplet $(j,m,\ell,k)$ is also
present, as well as the  other six obtained by circular
permutation (at least in the nondegenerate case mentioned in Lemma 
\ref{lem:rectangles}; adaptation to the degenerate case is obvious). 
These are the only occurrences of the corresponding
rectangle of wave numbers, and they produce the sum  
$$2 \( \overline a_ja_k\overline a_\ell a_m 
+ \overline a_ka_\ell\overline a_m a_j 
+ \overline a_\ell a_m\overline a_ja_k 
+ \overline a_ma_j\overline a_ka_\ell \) 
= 8 \RE \( \overline a_ja_k\overline
  a_\ell a_m \) ,$$
which is real. We consequently infer
\begin{equation*}
  \d_t\sum_{j\in J} |a_j(t,x+ t\kappa_j)|^2=0,  
\end{equation*}
and thus also
\begin{equation*}
  \frac{d}{dt}\sum_{j\in J} \|a_j(t,\cdot+t \kappa_j)\|_{L^2}^2
    = \frac{d}{dt}\sum_{j\in J} \|a_j(t,\cdot)\|_{L^2}^2 = 0. 
\end{equation*} 
\end{proof}
Let us now turn to the possibility of creating new profiles by
nonlinear interactions (note however that the conservation law
\eqref{eq:masscon} gives a global constraint on this process). To
simplify the presentation, we assume $d=2$. The creation of new 
oscillations in the general case $d\ge 2$ then 
follows by completing elements in $\R^2$ with $(0,\ldots,0)\in
\R^{d-2}$ and analogously for the situation on $\T^d$. Consider the
geometry associated to Example~\ref{ex:phase0}: We thus have (on
$\T^d$ or $\R^d$) 
\begin{equation*}
  i\d_t a_0 = \l \sum_{(k,\ell,m)\in I_0} a_k \overline a_\ell a_m.
\end{equation*}
Recall that $(k,\ell,m)\in I_0$ if and only if 
\begin{align*}
  \k_k-\k_\ell+\k_m= 0 \quad ; \quad
|\k_k|^2-|\k_\ell|^2+|\k_m|^2= 0,
\end{align*}
which obviously implies $\k_k\cdot \k_ m=0$. Such a
possibility occurs in two cases:
\begin{itemize}
\item$\k_k=0$ or $\k_m=0$.
\item $(\k_k,\k_m)=(\k_1,\k_3)$ or $(\k_k,\k_m)=(\k_3,\k_1)$ and hence
  $\k_\ell=\k_2$.  
\end{itemize}
From these various cases, we infer
\begin{equation*}
  i\d_t a_0 = \l \(|a_0|^2 +2|a_1|^2 +|a_2|^2 +2|a_3|^2\)a_0
  +2a_1\overline a_2 a_3.
\end{equation*}
Consider three non-vanishing initial oscillations, such that
$a_1\overline 
a_2 a_{3\mid t=0} \not =0$. Thus, even if $a_{0\mid t=0}=0$, we have
$\d_t a_{0\mid t=0}\not =0$, and this (non-oscillating) fourth mode is
instantaneously non-vanishing.  

\section{Higher order nonlinearities}
\label{sec:higher}

\subsection{Analysis of possible resonances}
\label{sec:reshigher}

So far we were only concerned with four-wave interactions
corresponding to cubic nonlinearities, \emph{i.e.} $\si =1$
in \eqref{eq:nls}. In general though, the set of resonances associated
to a (gauge invariant) nonlinearity of the form $f(z)=\l |z|^{2\si}z$,
$\si\in\N$, are defined by 
\begin{equation*}
  I_j^\si=\Big \{ \(\ell_1,\ldots,\ell_{2\si+1}\)\in J^{2\si+1};\
    \sum_{k=1}^{2\si+1}(-1)^{k+1} \k_{\ell_k} = \k_j,\  \sum_{k=1}^{2\si+1}(-1)^{k+1}
    |\k_{\ell_k}|^2 = |\k_j|^2\Big\}. 
\end{equation*}
As in Section \ref{sec:phase}, the set of wave vectors $\{\k_j\}_{j\in
  J}$ is constructed 
by induction, starting from an a finite or countable set
$\{\k_j\}_{j\in J_0}$, to  
which we first add a vector $\k$ when there exist
$\k_{\ell_1},\ldots,\k_{\ell_{2\si+1}} \in J_0$ such that  
\begin{equation} \label{eq:sires}
\sum_{k=1}^{2\si+1}(-1)^{k+1} |\k_{\ell_k}|^2 = 
\left|\sum_{k=1}^{2\si+1}(-1)^{k+1} \k_{\ell_k}\right|^2;
\end{equation}
we then set $\k = \sum_{k=1}^{2\si+1}(-1)^{k+1} \k_{\ell_k}$. The same iterative 
procedure as in \S\ref{sec:phasegen} leads to the following analogue
to Proposition  
\ref{prop:J}:
\begin{proposition} 
Let $\si\ge2$, and consider any $(2\si+1)$-tuple $(\k_{\ell_1},\ldots,
\k_{\ell_{2\si+1}})$ of wave vectors from $S=\{\kappa_j\}_{j\in J}$. Then, either 
the relation \eqref{eq:sires} is not satisfied, or the vector
$\k_j=\sum_{k=1}^{2\si+1}(-1)^{k+1} \k_{\ell_k}$ belongs to $S$.   
\end{proposition}

\begin{remark}
It is worth noting that, even if we only have very poor information on the set 
of wave vectors $\{\kappa_j\}_{j\in J}$, it is however a subset of the group 
generated by the initial set $\{\kappa_j\}_{j\in J_0}$.
\end{remark}

The profile equations, analogue to \eqref{eq:transportsystem}, are
then, for all $j\in J$:
\begin{equation}\label{eq:transportsystemsi}
 \d_t a_j +\k_j \cdot \nabla a_j 
  = -i \l\sum_{(\ell_1,\ldots,\ell_{2\si+1})\in I_j}a_{\ell_1} \overline
  a_{\ell_2} \dots a_{\ell_{2\si+1}}\quad ; \quad a_j(0,x) =    \alpha_j ( x).
\end{equation}

\subsection{The case of two modes}

Similar to the situation for $\si = 1$, the case of only two
initial modes is rather special. 
Indeed, the fact that two phases cannot create a new one 
extends also to higher order nonlinearities. In order 
to explain the argument, consider first a quintic nonlinearity,
corresponding to $\si = 2$.  
To obtain a nonlinear resonance, the wave vectors need to satisfy
\begin{align*}
  \k_k-\k_\ell+\k_m-\k_p+\k_q&=\k_j,\\
|\k_k|^2-|\k_\ell|^2+|\k_m|^2-|\k_p|^2+|\k_q|^2&=|\k_j|^2,
\end{align*}
where $k,\ell,m,p,q\in \{ j_1,j_2\}$, $j_1,j_2\in J$. First, if
$j_1$ (or $j_2$) 
appears at least twice on the left hand side, with at least one plus
and one minus, then the cancellation reduces the discussion to the one
we had about the cubic nonlinearity. Hence, no new resonant phase can be
created in this case. The complementary case corresponds, up to
exchanging $j_1$ and $j_2$, to
\begin{equation*}
  \k_k= \k_m= \k_q=\k_{j_1}\quad \text{and}\quad \k_\ell=\k_p= \k_{j_2}.
\end{equation*}
The above relations yield
\begin{equation*}
  3\k_{j_1}-2\k_{j_2}=\k_j\quad ;\quad 3|\k_{j_1}|^2-2|\k_{j_2}|^2=|\k_j|^2.
\end{equation*}
Squaring the first identity and comparing with the second one, we
infer
\begin{equation*}
  6|\k_{j_1}-\k_{j_2}|^2=0.
\end{equation*}
Therefore, no new resonant phase can be created by the quintic
interaction of two initial resonant plane waves. 

Consider now the general case where $\si
\ge 2$: The same argument as above shows that the only new case is the
one where all the plus signs correspond 
to one phase, and all the minus signs to the other:
\begin{equation*}
  (\si+1)\k_{j_1}-\si \k_{j_2}=\k_j\quad ;\quad
  (\si+1)|\k_{j_1}|^2-\si|\k_{j_2}|^2 =|\k_j|^2.  
\end{equation*}
Squaring the first identity and comparing with the second one, we
infer
\begin{equation*}
  \si(\si+1)|\k_{j_1}-\k_{j_2}|^2=0.
\end{equation*}
We conclude as above, and obtain the following result: 
%\begin{proposition}\label{prop:2modes}
%  Let $\si\in\N^*$, $\l\in \R^*$, and consider two phases 
%  \begin{equation*}
%    \phi_j(t,x)=\k_j\cdot x-\frac{t}{2}|\k_j|^2,\quad j=1,2,\ \k_1\not
%    =\k_2. 
%  \end{equation*}
%With $J=\{j_1,j_2\}$, we have, for  $k=1$ or $2$:
%\begin{align*}
%  I_{j_k}^\si = \big\{ &\(\ell_1,\ldots,\ell_{2\si+1}\)\in J^{2\si+1}\ ;\
%  \exists m \text{ odd,  }\k_{\ell_m}=\k_{j_k},\ \forall n\in
%  \{1,\ldots,2\si+1\}\setminus\{m\}, \\
%&\exists n'\in
%\{1,\ldots,2\si+1\}\setminus\{m\},\ n+n'\text{ odd},
%\k_{\ell_n}=\k_{\ell_{n'}}  \big\},
%\end{align*}
%and $I_j^\si=\emptyset$ otherwise. 
%\end{proposition}

\begin{proposition} \label{prop:2modes}
  Let $\si\in\N^*$, and let $\k_1,\k_2\in\R^d$ be such that
  $\k_1\neq\k_2$. To these wave vectors are associated the
  characteristic phases  
  \begin{equation*}
    \phi_j(t,x)=\k_j\cdot x-\frac{t}{2}|\k_j|^2,\quad j=1,2. 
  \end{equation*}
  Then, these two phases can not create no new phase by
  $(2\si+1)$th-order interaction: the set
  \begin{align*}
    \Big\{ \kappa\in\R^d \mid \exists (\ell_1,\ldots,\ell_{2\si+1})&\in
  \{1,2\}^{2\si+1}, \quad \k=\sum_{k=1}^{2\si+1}(-1)^{k+1}\k_{\ell_k}\\
  &\text{and}\quad |\k|^2=\sum_{k=1}^{2\si+1}(-1)^{k+1}|\k_{\ell_k}|^2) \Big\}
  \end{align*}
  is reduced to $\{\k_1,\k_2\}$.
\end{proposition}

In view of Proposition~\ref{prop:2modes}, the system 
\eqref{eq:transportsystemsi} becomes a
system of two equations, which can be integrated explicitly, as in
\cite[Remark~3.1]{CCTper}:  
\begin{equation}\label{eq:profilehigher}
  \begin{aligned}
    &\d_t a_j +\k_j\cdot\nabla a_j = -i\l \sum_{n=0}^\si \(
    \begin{array}[c]{c}
      \si+1 \\
n
    \end{array}\)
\(
    \begin{array}[c]{c}
      \si \\
n
    \end{array}\)|a_j|^{2\si-2n}|a_\ell|^{2n}
a_j,\\
&\d_t a_\ell +\k_\ell \cdot\nabla a_\ell = -i\l \sum_{n=0}^\si \(
    \begin{array}[c]{c}
      \si+1 \\
n
    \end{array}\)
\(
    \begin{array}[c]{c}
      \si \\
n
    \end{array}\)|a_\ell|^{2\si-2n}|a_j|^{2n}a_\ell.
  \end{aligned}
\end{equation}
In the case of $\T^d$, we find for instance
\begin{equation}\label{eq:profilehighersol}
 \begin{aligned}
  a_j(t)&= \alpha_j \exp\(-i \l t\sum_{n=0}^\si \(
    \begin{array}[c]{c}
      \si+1 \\
n
    \end{array}\)
\(
    \begin{array}[c]{c}
      \si \\
n
    \end{array}\)|\alpha_j |^{2\si-2n}|\alpha_\ell |^{2\ell}\),\\
a_\ell(t)&= \alpha_\ell \exp\(-i\l t\sum_{n=0}^\si \(
    \begin{array}[c]{c}
      \si+1 \\
n
    \end{array}\)
\(
    \begin{array}[c]{c}
      \si \\
n
    \end{array}\)|\alpha_\ell|^{2\si-2n}|\alpha_j |^{2n}\).
     \end{aligned}
\end{equation}
In the case of $\R^d$, the formula is more intricate and we shall omit it. 
\smallskip

Apart from the two-phase situation, the results for of Section 
\ref{sec:phasegen} on resonances do not
carry over to the general case $\si \ge 2$ in any straightforward
manner. Even in space dimension $d=1$, the resonant sets
cease to be as simple for $\si\ge 2$, provided that one starts with at
least three modes. 

\begin{example} Consider the quintic case $\si=2$ in $d=2$ spatial
  dimensions. As we have seen above a resonance for such a quintic
  nonlinearity appears if and only if
\begin{align*}
  \k_k-\k_\ell+\k_m-\k_p+\k_q&=\k_j,\\
|\k_k|^2-|\k_\ell|^2+|\k_m|^2-|\k_p|^2+|\k_q|^2&=|\k_j|^2.
\end{align*}
We can pick for instance three initial phases of the form
$$\k_1=(-1,0) \quad ; \quad \k_2=(0,0) \quad ;  \quad \k_3=(2,0) .
$$
For $k=1$, $\ell=p=2$, $m=q=3$, we have a resonance, creating
$\k_4=(3,0)$, whereas in the case $\si=1$, no resonance occurs between
the phases with wave vectors $\k_1$, $\k_2$ and $\k_3$. This example
shows that the geometric characterization of four-wave resonances given in \S\ref{sec:3phases}
does not export to  the case of six-wave resonances: $\k_1$, $\k_2$,
$\k_3$ and $\k_4$ all belong to the line $x_2=0$. 
\end{example}
\begin{example} Consider the same example as above in $d=1$. \emph{i.e.} 
pick three initial phases of the form
$$\k_1=-1 \quad ; \quad \k_2=0 \quad ;  \quad \k_3=2 
$$
and create a resonance $\k_4=3$ for $k=1$, $\ell=p=2$, $m=q=3$. 
This is in sharp contrast to the case $\si=1$, where no new phases can be created in $d=1$. 
Moreover, a non-vanishing amplitude $a_4$ is effectively generated:
\begin{align*}
(\d_t + \k_4\d_x) a_4 = 
& -3i \l \( |a_1|^4 + |a_3|^4 + |a_3|^4 + 
4(|a_1|^2|a_2|^2 + |a_2|^2|a_3|^2 + |a_3|^2|a_1|^2) \) a_4 \\
& -i 6 \l a_1 \overline{a_2}^2 a_3^2 - 6 i \l \overline a_1 |a_2|^2 a_3.
\end{align*}
We see that we may have $a_{4\mid t=0}=0$, but
\begin{equation*}
  \(\d_t a_4\)_{\mid t=0}= \(-i 6 \l a_1 \overline{a_2}^2 a_3^2 - 6 i
  \l \overline a_1 |a_2|^2 a_3\)_{\mid t=0}\not =0,
\end{equation*}
showing the appearance of a non-trivial $a_4$ for $t>0$. 
\end{example}

Despite this lack of knowledge concerning the precise structure of possible resonances for higher order 
nonlinearities, we shall see that 
we are able to prove the validity of WKB approximation even in this case. 

\section{Analytical framework}
\label{sec:analytical}

We now present the analytical framework needed for the rigorous
justification of a multiphase WKB approximation. 

\subsection{Wiener algebras}
\label{sec:wiener}

On $\mathcal M = \T^d$, we consider the usual Wiener algebra of functions with
absolutely summable Fourier series:  

\begin{definition}[Wiener algebra on $\mathcal M = \T^d$]
Functions of the form 
  \begin{equation*}
    f(y)=\sum_{k\in\Z} b_k e^{i \k_k\cdot y}, \quad \text{with }
    \k_k\in\Z^d \text{ and }
    b_k\in\C, 
  \end{equation*}
belong to $W(\T^d)$ if and only if $(b_k)_{k\in \Z}\in
\ell^1(\Z)$. We denote
\begin{equation*}
  \|f\|_{W} = \sum_{k\in \Z} |b_k|. 
\end{equation*}
\end{definition}
In the sequel, when $x\in \T^d$, we consider initial data for
\eqref{eq:nls} which are 
of the form $f(x/\eps)$, with $f\in W(\T^d)$  and $\eps^{-1}\in
\N^*$. 
\begin{lemma} \label{lem:substitT}
Let $f$ belong to $W(\T^d)$. Then, for all $\eps>0$ such that
$\eps^{-1}\in\N^*$, we have  $f(\cdot/\eps)\in W(\T^d)$, and 
$$\|f(\cdot/\eps)\|_W = \|f\|_W.$$
\end{lemma}
\smallbreak

For $\mathcal M = \R^d$, the framework is a bit different. Define the Fourier transform by
  \begin{equation*}
  \F f(\xi)=\widehat
  f(\xi)=\frac{1}{(2\pi)^{d/2}}\int_{\R^d}f(x)e^{-ix\cdot \xi} dx.
\end{equation*}
With this normalization, we have $\F^{-1}f(x) = \F f(-x)$. Following
\cite{JMRWiener} and  \cite{MColinLannes}, we use on $\R^d$ two
different Wiener-type algebras: For the exact solution we use $W(\R^d)$, \emph{i.e.} the space 
of functions with Fourier transform in $L^1(\R^d)$, and for
the profiles, we use $\A(\R^d)$, the space of almost periodic
$W(\R^d)$-valued functions on $\R^d$, with absolutely summable Fourier
series. We also set $\A(\T^d)=W(\T^d)$,
equipped with the same norm. 
\begin{definition}[Wiener algebra on $\mathcal M = \R^d$]
 We define
 \begin{equation*}
   W(\R^d)= \left\{ f\in {\mathcal S}'(\R^d;\C),\ \
     \|f\|_W:=\|\widehat f\|_{L^1(\R^d)} <\infty\right\}. 
 \end{equation*}
Functions of the form 
  \begin{equation*}
    f(x,y)=\sum_{k\in \Z} b_k(x) e^{i\k_k\cdot y}, \quad \text{with }
   \k_k\in\R^d \text{ and }  b_k\in W(\R^d) , 
 \end{equation*}
belong to $\A(\R^d)$ if and only if 
\begin{equation*}
   \|f\|_{\A}  := \sum_{k\in\Z}\|b_k\|_W=\sum_{k\in\Z}\|\widehat
   b_k\|_{L^1(\R^d)}<\infty.  
\end{equation*}
\end{definition}
In the sequel, when $x\in \R^d$, we consider initial data for
\eqref{eq:nls} which are 
of the form $f(x,x/\eps)$, with $f\in \A(\R^d)$. Again, we have
\begin{lemma} \label{lem:substitR}
Let $f\in \A(\R^d)$ and $\eps>0$. Then
$f(\cdot,\cdot/\eps)\in W(\R^d)$ and 
$$\|f(\cdot,\cdot/\eps)\|_W \le \|f\|_{\A}.$$
\end{lemma}

\begin{proof}
We simply have, when $f(x,y)=\sum_{k\in \Z} b_k(x) e^{i\k_k\cdot y}$:
$$\|f(\cdot,\cdot/\eps)\|_W = \big \lVert \sum_{k\in
    \Z}\widehat{b_k}(\cdot-\k_k/\eps)\big \rVert_{L^1(\R^d)}  
\le \sum_{k\in \Z} \|\widehat{b_k}(\cdot-\k_k/\eps)\|_{L^1(\R^d)} 
= \sum_{k\in \Z} \|\widehat{b_k}\|_{L^1(\R^d)}.$$
The last term is, by definition, $\|f\|_{\A}$. 
\end{proof}
\bigbreak

Denote (in the periodic setting as well as in the Euclidean case) 
\begin{equation*}
  U^\eps(t)=e^{i \eps \frac{t}{2}\Delta}.
\end{equation*}
The following properties will be useful (see \cite{MColinLannes}, and also
\cite{JMRWiener,Lan01,BL02}).
\begin{lemma} \label{lem:propWiener}
Let $\mathcal M=\T^d$ or $\R^d$.\\
$1.$ $W(\mathcal M)$ is a Banach space, continuously embedded into
$L^\infty(\mathcal M)$.\\ 
$2.$  $W(\mathcal M)$ is an algebra, in the sense that the
  mapping $(f,g)\mapsto fg$ is continuous from $W(\mathcal M)^2$ to
  $W(\mathcal M)$, and moreover
  \begin{equation*}
    \forall f,g\in W(\mathcal M), \quad
    \|fg\|_W\le \|f\|_W\|g\|_W.
  \end{equation*}
$3.$ If $F : \C\rightarrow\C$ maps $u$ to a finite sum of terms of the form
$u^{p}\overline u^q$, $p,q\in \N$, then it extends to a map from
$W(\mathcal M)$ to itself which is 
uniformly Lipschitzean on bounded sets of $W(\mathcal M)$. \\
$4.$ For all $t\in \R$, $U^\eps(t)$ is unitary on $W(\mathcal M)$. 
\end{lemma}

\subsection{Action of the free Schr\"odinger group on $W(\mathcal M)$}

As it is standard for solutions to the equation
\begin{equation*}
  i\eps\d_t w^\eps +\frac{\eps^2}{2}\Delta w^\eps = F^\eps,
\end{equation*}
we will consider the corresponding Duhamel's formula
\begin{equation*}
  w^\eps(t,x)=U^\eps(t)w^\eps(0,x)-
  i\eps^{-1}\int_0^tU^\eps(t-\tau)F^\eps(\tau,x)d\tau. 
\end{equation*}

In view of this representation formula we first need to study the
action of the free Schr\"odinger group $U^\eps(t)$ on $W(\mathcal
M)$. 

\subsubsection{The case $\mathcal M=\T^d$}

The action of $U^\eps(t)$ on Fourier series on $\T^d$ is well
understood. For $\sum_{k\in \Z} b_k e^{i \k_k\cdot y} \in W(\T^d)$: 
\begin{equation}\label{eq:identity}
  U^\eps(t)\(\sum_{k\in \Z} b_k e^{i \k_k\cdot x/\eps}\) = \sum_{k\in
    \Z} b_k
   e^{i\k_k\cdot x/\eps-i|\k_k|^2 t/(2\eps)}.  
\end{equation}
In view of Duhamel's formula, we will use the following
\begin{lemma}\label{lem:duhamelT}
  Let $T>0$, $\om\in\Z$, $\k\in\Z^d$, and $b,\d_t b\in L^\infty([0,T])$.
Denote
  \begin{equation*}
   D^\eps(t,x):= \int_0^t U^\eps(t-\tau)\( b(\tau) e^{i\k\cdot x/\eps-i\om
      \tau/(2\eps)}\)d\tau .
  \end{equation*}
$1.$ We have $D^\eps\in C([0,T]\times\T^d)$ and
\begin{equation*}
    \lVert D^\eps\rVert_{L^\infty([0,T]\times \T^d)}\le \int_0^T |b(t)|dt.
\end{equation*}
$2.$ Assume $\om \not = |\k|^2$. 
Then there exists $C$ independent of $\k$, $\om$ and $b$ such that
  \begin{equation*}
    \lVert D^\eps\rVert_{L^\infty([0,T]\times \T^d)}\le
    \frac{C\eps}{\left\lvert |\k|^2-\om\right\rvert}\(
    \|b\|_{L^\infty([0,T])} +\|\d_t b\|_{L^\infty([0,T])}\) . 
  \end{equation*}
\end{lemma}

\begin{proof}
  In view of the identity \eqref{eq:identity}, we have
  \begin{align*}
    D^\eps(t,x)&= \int_0^t b(\tau) e^{i\k\cdot x/\eps-i\om
      \tau/(2\eps)} e^{-i|\k|^2(t-\tau)/(2\eps)}d\tau\\
&= e^{i\k\cdot x/\eps-i|\k|^2t/(2\eps)}\int_0^t b(\tau)
e^{i(|\k|^2-\om)\tau/{2\eps}}d\tau. 
  \end{align*}
The first point is straightforward. 
Integration by parts yields, since by assumption  $|\k|^2-\om \in
\Z\setminus\{0\}$:  
with $\phi(t,x) = \k\cdot x - |\k|^2t/2$, 
\begin{align*}
 D^\eps(t,x) = e^{i\phi(t,x)/\eps}
 \Big(&-\frac{2\eps i}{|\k|^2-\om}b(\tau)e^{i(|\k|^2-\om)\tau/{2\eps}}
 \Big|_0^t\\
&
 + \frac{2\eps i}{|\k|^2-\om}\int_0^t
 \d_tb(\tau)e^{i(|\k|^2-\om)\tau/{2\eps}} d\tau\Big).
\end{align*}
The lemma then follows easily. 
\end{proof}
\subsubsection{The case $\mathcal M=\R^d$}

The Euclidean counterpart of Lemma~\ref{lem:duhamelT} is a little bit
more delicate: 

\begin{lemma}\label{lem:duhamelR}
 Let $T>0$, $\om\in\R$, $\k\in\R^d$, and $b\in L^\infty([0,T];W(\R^d))$. Denote
  \begin{equation*}
   D^\eps(t,x):= \int_0^t U^\eps(t-\tau)\( b(\tau,x) e^{i\k\cdot x/\eps-i\om
      \tau/(2\eps)}\)d\tau .
  \end{equation*}
$1.$ We have $D^\eps\in C([0,T];W(\R^d))$ and 
\begin{equation*}
  \lVert D^\eps\rVert_{L^\infty([0,T];W)}\le \int_0^T \left\lVert
      b(t,\cdot)\right\rVert_{W}dt.
\end{equation*}
$2.$ Assume $\om \not =|\k|^2$, and $\d_t b ,\Delta b \in L^\infty([0,T];W)$.
%.  Then
%\begin{equation}\label{eq:cvnorate}
%    \lVert D^\eps\rVert_{L^\infty([0,T];W)}\Tend \eps 0 0. 
%  \end{equation} 
%If in addition $\Delta b\in L^\infty([0,T];W)$, 
Then we have the control 
  \begin{equation*}
    \lVert D^\eps\rVert_{L^\infty([0,T];W)}\le
    \frac{C\eps}{\left\lvert |\k|^2-\om\right\rvert} \( \left\lVert
      b\right\rVert_{L^\infty([0,T];W)} +  
 \left\lVert {\Delta b}\right\rVert_{L^\infty([0,T];W)}+
\left\lVert {\d_t b}\right\rVert_{L^\infty([0,T];W)}\),
  \end{equation*} 
where $C$ is independent of $\k$, $\om$ and $b$. 
\end{lemma}
\begin{proof}
  By the definition of $U^\eps(t)$, we have
  \begin{equation*}
    \widehat D^\eps(t,\xi)= \int_0^t e^{-i\eps(t-\tau)|\xi|^2/2}
    \, \widehat b\(\tau,\xi-\frac{\k}{\eps}\)e^{-i\om
      \tau/(2\eps)}d\tau. 
  \end{equation*}
Setting $\eta = \xi -\k/\eps$, we have
 \begin{align*}
    \widehat D^\eps(t,\xi)&= e^{-i\eps t|\eta+\k/\eps|^2/2}\int_0^t
    e^{i\eps\tau|\eta+\k/\eps|^2/2} \,
    \widehat b\(\tau,\eta\)e^{-i\om
      \tau/(2\eps)}d\tau\\
&=e^{-i\eps t|\eta+\k/\eps|^2/2}\int_0^t
    e^{i\tau\theta/2} \,
    \widehat b\(\tau,\eta\)d\tau ,
  \end{align*}
where we have denoted 
\begin{equation*}
  \theta = \eps \left\lvert \eta+\frac{\k}{\eps}\right\rvert
  -\frac{\om}{\eps} = \underbrace{\eps|\eta|^2 +2\k\cdot
    \eta}_{\theta_1} +\underbrace{\frac{|\k|^2-\om}{\eps}}_{\theta_2}. 
\end{equation*}
The first point of the lemma is straightforward. To prove the second
point, integrate by parts, by first integrating $e^{i\tau\theta_2/2}$: 
\begin{equation*}
  \widehat D^\eps(t,\xi) = -\frac{2i}{\theta_2} e^{i\tau\theta/2} 
    \widehat b\(\tau,\eta\)\Big|_0^t +\frac{2i}{\theta_2}\int_0^t
    e^{i\tau\theta/2} \( i\frac{\theta_1}{2}\widehat
    b\(\tau,\eta\)+
\widehat {\d_t b}\(\tau,\eta\)\)d\tau.
\end{equation*}
We infer, if $b,\d_tb,\Delta b\in L^\infty([0,T];W)$:
\begin{equation*}
  \sup_{t\in [0,T]}
\lVert \widehat D^\eps(t)\rVert_{L^1} \lesssim
\frac{1}{|\theta_2|} \(
\lVert \widehat b\rVert_{L^\infty([0,T];L^1)} + 
 \lVert \widehat{\Delta b}\rVert_{L^\infty([0,T];L^1)}+
\lVert \widehat{\d_t b}\rVert_{L^\infty([0,T];L^1)}\).
\end{equation*}
This yields the second point of the lemma. 
%The weaker asymptotics
%\eqref{eq:cvnorate}  follows from a density argument. 
\end{proof}

\subsection{Construction of the exact solution} 
As a preliminary step in establishing a WKB approximation we first
need to know that \eqref{eq:nls} is well posed on $W(\mathcal M)$. 

\label{sec:constru}
\begin{proposition}\label{prop:existsol}
 Consider the initial value problem 
\begin{equation}\label{eq:nlsivp}
  i\eps\d_t u^\eps +\frac{\eps^2}{2}\Delta u^\eps = \l \eps
  |u^\eps|^{2\si}u^\eps\quad ;\quad u^\eps(0,x)=u_0^\eps(x),
\end{equation}
where $\si\in \N^*$, $\l\in \R$, and $x\in \mathcal M$, with either
$\mathcal M=\R^d$, 
or $\mathcal M=\T^d$, in which case $\eps^{-1}\in \N^*$. If
$u_0^\eps\in W(\mathcal M)$, then there exists $T^\eps>0$ and a unique
solution $u^\eps \in 
C([0,T^\eps];W(\mathcal M))$ to \eqref{eq:nlsivp}. 
\end{proposition}

\begin{remark}
  At this stage, the dependence of $T^\eps$ upon $\eps$ is unknown. In
  particular, $T^\eps$ might go to zero as $\eps\to 0$. The proof below
  actually shows that if $u_0^\eps$ is uniformly bounded in
  $W(\mathcal M)$ for $\eps\in ]0,1]$, then $T^\eps>0$ can be chosen
  independent of $\eps$. This case includes initial
  data \eqref{eq:CIgen} which we consider for the WKB analysis. 
\end{remark}
\begin{proof}
  Duhamel's formulation of \eqref{eq:nlsivp} reads
  \begin{equation*}
    u^\eps(t) = U^\eps(t)u_0^\eps - i\l\int_0^t
    U^\eps(t-\tau)\(|u^\eps|^{2\si}u^\eps(\tau)\)d\tau. 
  \end{equation*}
Denote by $\Phi^\eps(u^\eps)(t)$ the right hand side in the above
formula. From Lemmae~\ref{lem:propWiener}, \ref{lem:duhamelT} and
\ref{lem:duhamelR}, we have:
\begin{equation*}
  \left\lVert \Phi^\eps(u^\eps)(t)\right\rVert_{W}\le \|u_0^\eps\|_W
  + |\l| \int_0^t \|u^\eps(\tau)\|_W^{2\si+1}d\tau,
\end{equation*}
and if $\|u^\eps\|_{L^\infty([0,T];W)}, \|\widetilde
u^\eps\|_{L^\infty([0,T];W)}\le R$, then there exists $C=C(R)$ such
that
\begin{equation*}
  \left\lVert \Phi^\eps(u^\eps)(t)- \Phi^\eps(\widetilde
    u^\eps)(t)\right\rVert_{W}\le C(R)\int_0^t\left\lVert
    u^\eps(\tau)- \widetilde  u^\eps(\tau)\right\rVert_{W}
d \tau,\quad \forall t\in [0,T]. 
\end{equation*}
A fixed point argument in 
\begin{equation*}
  \left\{ u\in C([0,T];W(\mathcal M)),\ \sup_{t\in [0,T]}\|u(t)\|_W\le 2
    \|u_0^\eps\|_W\right\} 
\end{equation*}
for $T=T^\eps>0$ sufficiently small
then yields Proposition~\ref{prop:existsol}. 
\end{proof}

\subsection{Construction of the profiles}
\label{sec:constrprof}

In order to justify our multiphase WKB analysis, we first need to
establish an existence theory for the system of profile equations. To this end, for all $\si\in\N^*$, we rewrite 
the system~\eqref{eq:transportsystemsi} in its integral form: 
\begin{equation} \label{eq:transportsystemint}
\forall j\in J, \quad a_j(t,x) = a_j(0,x-t \k_j)  
  -i\l  \int_0^t \mathcal{N}_\si(a,\ldots,a)(\tau,x+(\tau-t)\k_j) d\tau,
\end{equation}
where, for $a^{(1)}=(a^{(1)}_j)_{j\in J}, \ldots , a^{(2\si+1)}=
(a^{(2\si+1)}_j)_{j\in J}$, we define the nonlinear term 
$\mathcal{N}_\si$ by: 
\begin{equation*}
\forall j\in J,\quad \mathcal{N}_\si \(a^{(1)},\ldots,a^{(2\si+1)}\)
  = \sum_{(\ell_1,\ldots,\ell_{2\si+1})\in I_j} a^{(1)}_{\ell_1} 
  \overline{a}^{(2)}_{\ell_2} \dots \overline
  a^{(2\si)}_{\ell_{2\si}}a^{(2\si+1)}_{\ell_{2\si+1}} .
\end{equation*}
It is clearly linear with respect to its arguments with odd exponents, 
and anti-linear with respect to the others. We prove in Lemma
\ref{lem:profNL} below that it is in fact well  
defined and continuous on $E(\mathcal M)$, for $\mathcal M=\T^d$ or $\mathcal M = \R^d$:
\begin{definition}
Define  
$$E(\R^d) = \{ a=(a_j)_{j\in J} \mid (\widehat{a}_j)_{j\in J}\in
\ell^1(J;L^1(\R^d)) \},$$
equipped with the norm
$$\|a\|_{E(\R^d)} = \sum_{j\in J}\|\widehat{a}_j\|_{L^1}.$$
Set also
$E(\T^d) = \ell^1(J)$, equipped with the usual norm 
$$
\|a\|_{E(\T^d)} = \sum_{j\in J}|a_j|.$$
\end{definition}

Note that $E$ simply represents, \emph{via} an isometric
correspondence, the family of coefficients of functions in $\A$ (up to
the choice of the wave numbers $\k_j$ in the case of $\R^d$): 
\begin{equation*}
 f(x,y) = \sum_{j\in J}a_j(x)e^{i\k_j\cdot y} \in \A(\R^d)
\text{ iff } a \in E(\R^d), 
\end{equation*}
and then $ \| a \|_{E}=\|f\|_{\A} $. The same holds for $\mathcal M = \mathbb T^d$.
%\begin{equation*}
 %f(y) = \sum_{j\in J}a_je^{i\k_j\cdot y} \in \A(\T^d) \text{ iff } a
%\in E(\T^d),  \quad \text{with $\|a\|_{E}=\|f\|_{\A}$.}
%\end{equation*}

\begin{lemma}
\label{lem:profNL}
Let $\si\in\N^*$. For $\mathcal M=\R^d$ or $\mathcal M=\T^d$, the nonlinear 
expression $\mathcal{N}_\si$ defines a continuous mapping  
from $E(\mathcal M)^{2\si+1}$ to $E(\mathcal M)$, and for all
$a^{(1)},\ldots,a^{(2\si+1)}\in E(\mathcal M)$
\begin{equation*}
  \left\lVert \mathcal{N}_\si\(a^{(1)},\ldots,a^{(2\si+1)}\) \right\rVert_{E} \le
  \|a^{(1)}\|_{E} \dots  
\|a^{(2\si+1)}\|_{E}.
\end{equation*}
\end{lemma}
\begin{proof}
We consider the case $\mathcal M=\R^d$, since $\mathcal M=\T^d$ is even 
simpler. In order to bound  
\begin{align*}
\lA \mathcal{N}_\si\(a^{(1)},\ldots,a^{(2\si+1)}\) \rA_{E} =\\
 = \sum_{j\in J} \Bigg\lVert \sum_{(\ell_1,\ldots,\ell_{2\si+1})\in I_j} &
\F\(a^{(1)}_{\ell_1}\) \ast \F\(\overline a^{(2)}_{\ell_2}\) \ast \dots 
\ast \F\(a^{(2\si+1)}_{\ell_{2\si+1}}\)\Bigg\rVert_{L^1} \\
 \le  \sum_{j\in J} \sum_{(\ell_1,\ldots,\ell_{2\si+1})\in I_j} & 
\lA \F\(a^{(1)}_{\ell_1}\) \ast \F\(\overline a^{(2)}_{\ell_2}\) \ast \dots 
\ast \F\(a^{(2\si+1)}_{\ell_{2\si+1}}\)\rA_{L^1} ,
\end{align*}
we use Young's inequality and observe that,  
once $j$, $\ell_1$, \dots, $\ell_{2\si}$ are chosen, $\ell_{2\si+1}$
is determined  
(since $\k_{\ell_{2\si+1}}=\k_j-\sum_{k=1}^{2\si} (-1)^{k+1}\k_{\ell_k}$, and 
$n\neq m \Rightarrow \k_n\neq\k_m$), so that  
\begin{equation*}
\lA \mathcal{N}_\si\(a^{(1)},\ldots,a^{(2\si+1)}\) \rA_{E} 
\le \sum_{(\ell_1,\ldots,\ell_{2\si+1})\in J^{2\si+1}}
\lA \F\(a^{(1)}_{\ell_1}\)\rA_{L^1}  
\dots \lA \F\(a^{(2\si+1)}_{\ell_{2\si+1}}\)\rA_{L^1},  
\end{equation*}
which gives the desired result.
\end{proof}

This consequently yields the following existence result for
\eqref{eq:transportsystemint}, where here and in the following we denote $\<\k\>^2\equiv 1 +
|\k|^2 $. 

\begin{proposition} \label{prop:existprof}
Let $\sigma\in\N^*$, and $\mathcal M=\R^d$ or $\mathcal M=\T^d$.\\
For all
$\alpha=(\alpha_j)_{j\in J}\in E(\mathcal M)$, there exist $T>0$ and a
unique solution 
\begin{equation*}
  t \mapsto a(t)=(a_j(t))_{j\in J}\in
C([0,T],E(\mathcal M))
\end{equation*}
to the  
system~\eqref{eq:transportsystemint}, with $a(0)=\alpha$. 
Moreover, the following properties hold:\\
$1.$ If
$(\<\k_j\>^s\alpha_j)_{j\in J}\in E(\mathcal M)$  for some $s \in\N$,
then $(\<\k_j\>^s a_j)_{j\in J}\in C([0,T],E(\mathcal
M))$.\\
$2.$ On $\mathcal M=\R^d$, if
$(\<\k_j\>^s\d^\beta_x\alpha_j)_{j\in J}\in E(\R^d)$, for some
$\beta\in\N^d$ and  $s\in\N$, then $(\<\k_j\>^s\d^\beta_x
a_j)_{j\in J} \in C([0,T];E(\R^d))$.  
\end{proposition} 

\begin{proof}
The existence result follows from Lemma~\ref{lem:profNL} and the
standard Cauchy--Lipschitz result for ODE's. Concerning the propagation
of moments $\<\k_j\>^s a_j$, we again apply a fixed-point argument,
estimating nonlinear terms 
$$\<\k_j\>^s \mathcal{N}_\si(a^{(1)},\ldots,
a^{(2\si+1)})$$ 
as in the proof of Lemma~\ref{lem:profNL}, \emph{via} 
\begin{align*}
\<\k_j\>^2\equiv 1 + |\k_j|^2 
& = 1 + \sum_{k=1}^{2\si+1} (-1)^{k+1} |\k_{\ell_k}|^2 \\
& \le \sum_{k=1}^{2\si+1} \<\k_{\ell_k}\>^2  \le (2\si+1) \prod_{k=1}^{2\si+1} \<\k_{\ell_k}\>^2 , 
\end{align*}
when $(\ell_1,\dots,\ell_{2\si+1})\in I_j$. 
The last statement of the proposition is concerned with the smooth 
dependence upon the parameter $x$. This follows by commuting  
\eqref{eq:transportsystemsi} with $\d_x$ and using the fact that
$W(\R^d)$ is an algebra, continuously embedded in $L^\infty$, since
then  
$$
\frac{d}{dt} \| \d_x a\|_{E} \lesssim \| \d_x \alpha \|_{E} + C(\| a
\|_{E}) \| \d_x a \|_{E} ,
$$
and a Gronwall argument shows that $ \| \d_x a\|_{E}$ remains bounded
for all $t \in [0,T]$. Similarly we conclude for the higher order
derivatives, possibly multiplied by weights $\<\k_j\>^s$. 
\end{proof}

For the particular situation for $\si=1$, in $d=1$ and/or the case 
of only two initial phases, we infer a stronger result, thanks to the 
explicit formulas given in \S\ref{sec:1d} and \S\ref{sec:2profiles}. 

\begin{corollary}
Under the assumption of Proposition~\ref{prop:existprof}, in the case 
$\si=1$, if in addition $d=1$, then $T$ can be taken arbitrarily large, 
with $a_j(t)$ explicitly given by 
\eqref{eq:selfmodulation} and \eqref{eq:selfmodulationformula}. 
Similarly, if $\sharp J_0 \le 2$, then $T$ can be taken arbitrarily large.
\end{corollary}
\begin{remark}
  In the case of higher order nonlinearities, 
  \emph{i.e.} $\si\ge 2$, Equation~\eqref{eq:profilehigher} makes it
  possible to see, 
  \emph{via} explicit integration (see \eqref{eq:profilehighersol} in
  the case of the torus), that if $\alpha_j,\alpha_\ell\in W(\mathcal
  M)$, then $a_j,a_\ell\in C([0,\infty[,W(\mathcal M))$. 
\end{remark}

\section{Rigorous justification of the multiphase WKB analysis}
\label{sec:justif}

\subsection{Construction of an approximate solution}

We start from oscillating initial data, given by a profile in 
$\A(\mathcal M)$, with $\mathcal M=\T^d$
or $\R^d$: 
$$u^\eps_{\rm app}(0,x)= \sum_{j\in
  J_0}\alpha_j(x) e^{i \kappa_j \cdot x/\eps},$$
with $\alpha_j(x)={\rm Const.}$ in the case $\mathcal M=\T^d$. 
\begin{assumption} \label{as:regularity}
For both $\mathcal M = \R^d$ and $\mathcal M= \mathbb T^d$ we assume $(\alpha_j)_{j\in J_0}\in E(\mathcal M)$. 
For $\mathcal M = \R^d$ we assume in addition
\begin{equation*} \forall |\beta|\le 2, \
(\partial^\beta_x\alpha_j)_{j\in J_0}\in E(\R^d),  \quad \text{and} \quad
\forall |\beta|\le 1, \ (\<\k_j\>\partial^\beta_x\alpha_j)_{j\in J_0}
\in E(\R^d).
\end{equation*}
\end{assumption}
From Proposition \ref{prop:existprof} we know, that these data produce a solution
$(a_j)_{j\in J}\in C([0,T],E(\mathcal M))$ to the amplitude system and we consequently 
define the approximate solution $u_{\rm app}^\eps$ by
\begin{equation}
  \label{eq:uapp}
  u_{\rm app}^\eps(t,x) = \sum_{j\in J}a^\eps_j(t,x) e^{i\phi_j(t,x)/\eps},
\end{equation}
with $\phi_j$ given by \eqref{eq:planewave}. 
The sequence $(a_j)_{j\in J}$ is such that
\begin{align*}
  &\(\partial^\beta_xa_j\)_{j\in J}\in C([0,T],E(\mathcal M)),
  \quad |\beta|\le 2,\\
& \(\<\k_j\> \partial^\beta_xa_j\)_{j\in J}\in
C([0,T],E(\mathcal M)), \quad |\beta|\le 1.
\end{align*}
We see from equation \eqref{eq:transportsystemint}  that 
$(\partial_ta_j)_{j\in J}\in C([0,T],E(\mathcal M))$. 
We find (in the sense of distributions)
\begin{equation*}
  i\eps \d_t u_{\rm app}^\eps +\frac{\eps^2}{2}\Delta u_{\rm app}^\eps 
  = \l\eps |u_{\rm app}^\eps|^{2\si} u_{\rm app}^\eps  - \l \eps r_1^\eps 
  + \eps^2 r_2^\eps,
\end{equation*}
where 
\begin{equation}\label{eq:reminder2}
  r_2^\eps=   \frac{1}{2} \sum_{j\in J} e^{i\phi_j/\eps}  \Delta a_j  ,
\end{equation}
and the remainder $r_1^\eps$ takes into account the non-characteristic
phases created by nonlinear interaction. This means that it is a sum
of terms of the form 
$$a_{\ell_1}\overline{a}_{\ell_2}\dots \overline a_{\ell_{2\si}} a_{\ell_{2\si+1}}
e^{i(\phi_{\ell_1}-\phi_{\ell_2}+\, \dots\, -\phi_{\ell_{2\si}}+\phi_{\ell_{2\si+1}})/\eps},$$ 
where the rapid phase is given by
\begin{equation*}
  \sum_{p=1}^{2\si+1}(-1)^{p+1}\phi_{\ell_p}(t,x) =
  \(\sum_{p=1}^{2\si+1}(-1)^{p+1}\k_{\ell_p}\)\cdot x -\frac{t}{2}
  \sum_{p=1}^{2\si+1}(-1)^{p+1}\lvert \k_{\ell_p} \rvert^2,
\end{equation*}
and
\begin{equation*}
  \left\lvert \sum_{p=1}^{2\si+1}(-1)^{p+1}\k_{\ell_p}\right\rvert^2
  \neq  \sum_{p=1}^{2\si+1}(-1)^{p+1}\left\lvert
    \k_{\ell_p}\right\rvert^2. 
\end{equation*}
In other words,  $(\ell_1,\dots,\ell_{2\si+1})$ belongs to the 
\emph{non-resonant set}
$$N := J^{2\si+1} \setminus \bigcup_{j\in J} I^\si_j.$$
With these conventions, we have 
\begin{equation}\label{eq:reminder1}
 r_1^\eps = \sum_{(\ell_1,\dots,\ell_{2\si+1})\in N}
 a_{\ell_1}\overline{a}_{\ell_2}\dots a_{\ell_{2\si+1}}  
 e^{i(\phi_{\ell_1}-\phi_{\ell_2}+ \, \dots\, -\phi_{\ell_{2\si}}+\phi_{\ell_{2\si+1}})/\eps} . 
\end{equation}
Estimating $r_2^\eps$ in $W$ is straightforward, since
$(\partial^\beta_xa_j)_{j\in J}\in C([0,T],E)$ for
$|\beta|\le 2$: 
\begin{equation} \label{eq:estimr2}
\|r^\eps_2\|_{W} \le \frac{1}{2} \|\Delta a \|_{E}.
\end{equation}  
Note that $r_2$ simply vanishes if $\mathcal M=\T^d$. 
In order to estimate $r_1^\eps$, we impose the following condition on
the set of wave numbers $\{\k_j\}_{j\in J}$. 
\begin{assumption}
\label{as:nosmalldiv}
There exists $c>0$ such that for all $(\ell_1,\dots,\ell_{2\si+1}) 
\in N $,  
\begin{equation*}
  \delta\(\ell_1,\dots,\ell_{2\si+1}\) \equiv
\left\lvert \left\lvert
    \sum_{p=1}^{2\si+1}(-1)^{p+1}\k_{\ell_p}\right\rvert^2
  -\sum_{p=1}^{2\si+1}(-1)^{p+1}\left\lvert \k_{\ell_p}\right\rvert^2
\right\rvert  \ge c.
\end{equation*}
\end{assumption}
\begin{remark}
$(i)$ This assumption is of course satisfied when only finitely many phases are created $\sharp J <\infty$. \\
$(ii)$ Similarly, this assumption holds for $\{\k_j\}_{j\in J} \subset
\Z^d$, since in this case, the quantity considered is an integer.  \\
$(iii)$ Consider the cubic case $\si=1$, and suppose that
$\{\k_j\}_{j\in J}$ is included in a rectangular 
net. Up to translation, this  rectangular net has the form  
$$\{Am \in \R^d \mid m \in \Z^d\},$$
with $A$ a $d\times d$ matrix of the form $A=RD$, where $D$ is
diagonal, and $R$ is a rotation. Then we have, for all $k, l, m \in
\Z^d$:  
\begin{align*}
\la \la Ak-A\ell+Am\ra^2-|Ak|^2+|A\ell|^2-|Am|^2\ra  & =
\la(Ak-A\ell)\cdot(Ak-Am)\ra\\ 
&  =\la(k-l)\cdot\((A^TA)(k-m)\)\ra.
\end{align*}
Since $^TAA=D^2$, denoting $\mu_1^2,\dots,\mu_d^2$ the squares of the
eigenvalues of $D$, Assumption~\ref{as:nosmalldiv} is then satisfied
if and only if the group generated by $\mu_1^2,\dots,\mu_d^2$ in $\R$ is
discrete, \emph{i.e.} these numbers are (pairwise) rationally
dependent.  
\end{remark}
The reason for imposing the above assumption is a small divisor
problem, as will become clear from the proof of the following
lemma. It is possible to relax Assumption~\ref{as:nosmalldiv} to
a less rigid one, to the cost of a more technical presentation. The
latter is sketched in an appendix.

\begin{lemma}\label{lem:Rest}
For $\mathcal M=\T^d$ or $\mathcal M=\R^d$, let $r^\eps_1$ be given by \eqref{eq:reminder1} and denote 
  \begin{equation}\label{eq:R}
   R_1^\eps(t,x):= \int_0^t U^\eps(t-\tau) r_1^\eps(\tau,x) d\tau, \quad \text {on $[0,T]\times \mathcal M$} .
  \end{equation}
Let Assumptions \ref{as:regularity}--\ref{as:nosmalldiv} hold. Then, there exists a constant
$C>0$, such that: 
\begin{equation*}
    \lVert R_1^\eps\rVert_{L^\infty([0,T];W(\mathcal M))}\le C \eps  .
  \end{equation*}
\end{lemma}
\begin{proof} 
We only treat the case on $\mathcal M = \R^d$ in detail. The case
$\mathcal M = \T^d$ can be treated analogously. We have  
\begin{align*}
   &R_1^\eps
    (t,x) = \\
   & \sum_{(\ell_1,\dots,\ell_{2\si+1})\in N}   
   \int_0^t U^\eps(t-\tau)
   \left( 
   (a_{\ell_1}\overline{a}_{\ell_2}\dots a_{\ell_{2\si+1}}) 
   e^{i\(\phi_{\ell_1}-\phi_{\ell_2}+\dots+\phi_{\ell_{2\si+1}}\)  
   /\eps} \right)(\tau,x) \, d\tau . 
\end{align*}
Thus, setting $b_{\ell_1,\dots,\ell_{2\si+1}}: = 
a_{\ell_1}\overline{a}_{\ell_2}\dots a_{\ell_{2\si+1}}$,
Lemma~\ref{lem:duhamelR} yields
\begin{align*}
 \lVert R_1^\eps\rVert_{L^\infty([0,T] ; W )}\lesssim \, & \, \eps
 \sum_{(\ell_1,\dots,\ell_{2\si+1})\in N} 
    \frac{1}{\delta(\ell_1,\dots,\ell_{2\si+1})}
    \Big( \| \widehat b_{\ell_1,\dots,\ell_{2\si+1}}\|_{L^\infty([0,T]; L^1)} \\ 
    & \hspace{0.7cm} + \| \widehat {\Delta  b}_{\ell_1,\dots,\ell_{2\si+1}}
    \|_{L^\infty([0,T]; L^1)} + 
    \| \widehat {\d_t
      b}_{\ell_1,\dots,\ell_{2\si+1}}\|_{L^\infty([0,T]; L^1)}\Big )
    \\  
    \lesssim  \, & \, \eps \sum_{(\ell_1,\dots,\ell_{2\si+1})\in N}
    \Big( \| \widehat b_{\ell_1,\dots,\ell_{2\si+1}}\|_{L^\infty([0,T]; L^1)} \\ 
    & \hspace{0.7cm} + \| \widehat {\Delta  b}_{\ell_1,\dots,\ell_{2\si+1}}
    \|_{L^\infty([0,T]; L^1)} + \| \widehat {\d_t b}_{\ell_1,\dots,\ell_{2\si+1}}
    \|_{L^\infty([0,T]; L^1)}\Big ),
\end{align*}
where we have used  Assumption~\ref{as:nosmalldiv}. 
Next, using Young's inequality, as in the proof of
Lemma~\ref{lem:profNL}, we get:  
\begin{align*}
\sum_{(\ell_1,\dots,\ell_{2\si+1})\in N} 
\| \widehat{b}_{\ell_1,\dots,\ell_{2\si+1}}\|_{L^\infty([0,T]; L^1)}
& \lesssim  \sum_{(\ell_1,\dots,\ell_{2\si+1})\in N} 
 \| \widehat{ a}_{\ell_1}\|_{L^\infty_T L^1} \dots 
 \| \widehat{ a}_{\ell_{2\si+1}}\|_{L^\infty_TL^1} \\
 \lesssim &\sum_{(\ell_1,\dots,\ell_{2\si+1})\in J^{2\si+1}} 
 \| \widehat{ a}_{\ell_1}\|_{L^\infty_T L^1} \dots 
 \| \widehat{ a}_{\ell_{2\si+1}}\|_{L^\infty_T L^1} \\
& \lesssim \|( a_j)_{j\in J}\|^{2\si+1}_{L^\infty([0,T];E)} .
\end{align*}
Leibniz formula and H\"older inequality yield similar
estimates for 
$\widehat {\Delta  b}_{\ell_1,\dots,\ell_{2\si+1}}$ 
and $\widehat {\d_t b}_{\ell_1,\dots,\ell_{2\si+1}}$ in
$L^\infty([0,T]; L^1(\mathcal M))$, and the lemma follows.
\end{proof}

\subsection{Accuracy of the multiphase WKB approximation}

With the above results in hand, we can now prove our main theorem.

\begin{theorem}[General approximation result]
\label{theo:justif}
Let $\si\ge1$, $\mathcal M=\T^d$ or $\R^d$, and Assumptions 
\ref{as:regularity}--\ref{as:nosmalldiv}  hold. Given an approximate solution $u^\eps_{\rm app}\in
C([0,T];W(\mathcal M))$ as in 
\eqref{eq:uapp}, we consider a family of initial data $(u^\eps_0)_{\eps>0} \in 
W(\mathcal M)$, such that   
$$ 
\left\lVert u^\eps_0-{u_{\rm app\mid t=0}^\eps}\right\rVert_{W(
  \mathcal M)}\le C_0 \eps, 
$$ 
for some $C_0\ge 0$ independent of $\eps$. 
Then there exists $\eps_0(T)>0$, such that for any $0<\eps\le \eps_0(T)$, the exact
solution to the Cauchy problem \eqref{eq:nlsivp} satisfies
$u^\eps \in L^\infty([0,T];W(\mathcal M))$. In addition, $u_{\rm
  app}^\eps$ approximates $u^\eps$ up  
to $\O(\eps)$:
\begin{equation*}
    \left\lVert u^\eps-u_{\rm
        app}^\eps\right\rVert_{L^\infty([0,T]\times \mathcal M)}\le  
\left\lVert u^\eps-u_{\rm
    app}^\eps\right\rVert_{L^\infty([0,T];W(\mathcal M))}\le C\eps, 
  \end{equation*}
where $C$ is independent of $\eps$. 
\end{theorem}
Obviously the result for $x \in \T$, announced in the introduction,
can be seen as a special case of Theorem \ref{theo:justif}. 

\begin{proof} 
From Proposition~\ref{prop:existsol}, we may consider a solution 
$u^\eps\in C([0,T^\eps],W(\mathcal M))$ to \eqref{eq:nls}. 
We define the difference $w^\eps:= u^\eps-u_{\rm app}^\eps$. Then 
$w^\eps \in C([0,\tau^\eps],W(\mathcal M))$, where $\tau^\eps=\min
(T^\eps,T)$. We prove that for $\eps$ sufficiently small, $w^\eps$ 
may be extended up to time $T$, with $w^\eps \in C([0,T],
W(\mathcal M))$. Take $\eps_0>0$ so that $C_0\eps_0\leq1/2$, and for 
$\eps\in]0,\eps_0]$, let  
$$t^\eps := \sup \Big\{ t\in[0,T] \mid
  \sup_{t'\in[0,t]}\|w^\eps(t')\|_{W(\mathcal M)}\le 1\Big\}.$$
We already know that $t^\eps>0$ by the local existence result for
$u^\eps$. By possibly reducing $\eps_0>0$, we shall show that $t^\eps\ge T$. The
error term $w^\eps$ solves: 
\begin{equation*}
  i \d_t w^\eps +\frac{\eps}{2}\Delta w^\eps = \l
 \(|u_{\rm app}^\eps+w^\eps|^{2\si} (u_{\rm app}^\eps+w^\eps) 
 - |u_{\rm app}^\eps|^{2\si} u_{\rm app}^\eps\) + \l  r_1^\eps - \eps r^\eps_2, 
\end{equation*}
where $r_1^\eps$, $r^\eps_2$ are given in
\eqref{eq:reminder2}--\eqref{eq:reminder1}.  
Using Duhamel's formula we can rewrite this equation as 
\begin{align*}
    w^\eps(t) = & \, U^\eps(t)w_0^\eps - i\l\int_0^t
    U^\eps(t-\tau)\(|u_{\rm app}^\eps+w^\eps|^{2\si} 
    (u_{\rm app}^\eps+w^\eps) -  |u_{\rm app}^\eps|^{2\si} 
    u_{\rm app}^\eps\)(\tau) \, d\tau \\ 
   & \, - i \l R_1^\eps(t) + i \eps \int_0^t U^\eps(t-\tau)
   r^\eps_2(\tau) \, d \tau ,  
  \end{align*}
  where $R_1^\eps$ is defined in \eqref{eq:R}. Using the fact that
  $U^\eps(t)$ is unitary on $W(\mathcal M)$, and the estimates given in
  \eqref{eq:estimr2} and in  Lemma~\ref{lem:Rest}, we obtain on $[0,t^\eps]$:
  \begin{align*}
  \lVert w^\eps(t)\rVert_{W(\mathcal M )} \le & \ C_1 \eps  + |\l|
  \int_0^t  \lVert  \(|u_{\rm app}^\eps+w^\eps|^{2\si} 
  (u_{\rm app}^\eps+w^\eps) -  |u_{\rm app}^\eps|^{2\si} 
  u_{\rm app}^\eps\)(\tau)
  \lVert_{W(\mathcal M )} \, d\tau\\ 
   \le & \ C_1 \eps +  C_2 \int_0^t  \lVert  w^\eps (\tau)
   \lVert_{W(\mathcal M )} \, d\tau, 
    \end{align*}
by the Lipschitz property from Lemma \ref{lem:propWiener} . Note that,
in view of Lemma~\ref{lem:substitT} , resp. Lemma~\ref{lem:substitR},
$(u^\eps_{\rm app})_{\eps>0}$ is a bounded family in
$C([0,T],W(\mathcal M))$, and restricting $t$ to
$[0,t^\eps]$ ensures that $w^\eps(t)$ stays bounded in $W(\mathcal
M)$. The constants $C_1$, $C_2$ depend on $C_0$ and $u^\eps_{\rm
  app}$. Now, Gronwall lemma yields  
$$\lVert w^\eps(t)\rVert_{W(\mathcal M )} \le C_1\eps \( 1 +
\frac{e^{C_2T}}{C_2} \),$$ 
and we may reduce $\eps_0$ so that $C_1\eps_0 \( 1 + e^{C_2T}/C_2 \) <
1$. This shows that $t^\eps\ge T$, for all $\eps\in ]0,\eps_0]$. Then,
$T^\eps\ge T$ follows, as well as the desired approximation of $u^\eps$
by $u^\eps_{\rm app}$, since $w^\eps=\O(\eps)$ in $L^\infty([0,T];W)$.      
\end{proof}

\section{Proof of the instability result}
\label{sec:instab}

This section is devoted to the proof of Theorem~\ref{theo:instab}. To this end  
we essentially rewrite the proof of M.~Christ,
J.~Colliander, and T.~Tao \cite{CCTper} in terms of weakly nonlinear
geometric optics. It then becomes easy to see that the justification
given in the previous paragraph makes it possible to extend the
one-dimensional analysis of \cite{CCTper} in order to infer
Theorem~\ref{theo:instab}. 

\begin{proof}[Proof of Theorem~\ref{theo:instab}]
We start with two Fourier modes, one of them being zero:
\begin{equation*}
  i\d_t u +\frac{1}{2}\Delta u =\l|u|^{2\si}u\quad ;\quad
  u(0,x)= \alpha_0 +\alpha_1 e^{iKx_1},\quad K\in \N. 
\end{equation*}
The fact that we privilege oscillations with respect to the first space
variable is purely arbitrary. Define $\widetilde u$ as the solution to
the same equation, with data
\begin{equation*}
  \widetilde u(0,x)= \widetilde \alpha_0 +\widetilde \alpha_1 e^{iKx_1}.
\end{equation*}
Let
\begin{equation*}
  \eps = \frac{1}{K^2}\quad ;\quad u^\eps(t,x)= u\(
  t,\frac{x}{\sqrt \eps}\) = u\(
  t,Kx\).
\end{equation*}
($\eps$ is chosen so that we remain on the torus.)
We see that $u^\eps$ solves \eqref{eq:nls} on $\T^d$, with
\begin{equation*}
  u^\eps(0,x)=\alpha_0 +\alpha_1 e^{ix_1/\eps}. 
\end{equation*}
From Theorem~\ref{theo:justif}, we know that there exists $T>0$ 
independent of $\eps$, such that
\begin{equation*}
  \|u^\eps - u_{\rm app}^\eps\|_{L^\infty([0,T]\times\T^d)}
+\|\widetilde u^\eps - \widetilde
u_{\rm app}^\eps\|_{L^\infty([0,T]\times\T^d)}=\O(\eps), 
\end{equation*}
where $u_{\rm app}^\eps$ is the approximate solution defined by \eqref{eq:uapp},
and $\widetilde u_{\rm app}^\eps$ is defined similarly. On the other hand, we
have
\begin{equation*}
  u_{\rm app}^\eps(t,x)= \alpha_0 e^{-i\l t\theta_0} + \alpha_1
  e^{-i\l t\theta_1}e^{i(x_1-t/2)/\eps}, 
\end{equation*}
where, in view of \eqref{eq:profilehigher}, $\theta_0$ is given
by
\begin{equation*}
  \theta_0 = \sum_{n=0}^\si \(
    \begin{array}[c]{c}
      \si+1 \\
n
    \end{array}\)
\(
    \begin{array}[c]{c}
      \si \\
n
    \end{array}\)|\alpha_0|^{2\si-2nl}|\alpha_1|^{2n}.
\end{equation*}
We infer, uniformly in
$t\in [0,T]$,
\begin{equation*}
  \left\lvert \int_{\T^d}(u(t,x)-\widetilde u(t,x))dx\right\rvert
  =\left\lvert \alpha_0 e^{-i\l t \theta_0} - \widetilde \alpha_0
    e^{-i\l t
      \widetilde\theta_0}\right\rvert 
  +\O(\eps),
\end{equation*}
with obvious notations. 

To prove the first point of Theorem~\ref{theo:instab}, set 
\begin{equation*}
  \alpha_0 =\widetilde \alpha_0=\frac{\rho}{2},\quad \alpha_1 =
  \frac{\rho}{2K^s} 
  =\frac{\rho}{2} K^{|s|},\quad  \widetilde \alpha_1 = \sqrt{\alpha_1^2
  +\frac{1}{\delta}}. 
\end{equation*}
We infer, for $0<\delta\le 1$,
\begin{equation*}
  \left\lvert \theta_0 - 
      \widetilde\theta_0\right\rvert \gtrsim \frac{1}{\delta}. 
\end{equation*}
We have $\|u(0)-\widetilde u(0)\|_{H^s}<\delta$ 
provided $K> \delta^{1/s}$. 
Since $\widetilde \alpha_0=\alpha_0$, we also have
\begin{equation*}
  \left\lvert \int_{\T^d}(u(t,x)-\widetilde u(t,x))dx\right\rvert
  =\left\lvert 2\alpha_0\sin
    \(\frac{\l t}{2}\(\widetilde\theta_0-\theta_0\)\) \right\rvert 
  +\O(\eps).
\end{equation*}
We infer that we can find $t\in [0,\delta]$ so that the right hand
side is bounded from below by $\rho/2$, provided $N$ is sufficiently
large (hence $\eps$ sufficiently small). 
\smallbreak

To prove the second point of Theorem~\ref{theo:instab}, set 
\begin{equation*}
  \alpha_0 =\frac{\rho}{2},\quad \widetilde \alpha_0=
  \alpha_0+\delta,\quad \alpha_1 =\widetilde \alpha_1 = 
  \frac{\rho}{2K^s} 
  =\frac{\rho}{2} K^{|s|}. 
\end{equation*}
For $\delta$ small compared to $\rho$, we use the same estimate as
above,
\begin{equation*}
  \left\lvert \int_{\T^d}(u(t,x)-\widetilde u(t,x))dx\right\rvert
  \gtrsim \left\lvert 2\alpha_0\sin
    \(\frac{\l t}{2}\(\widetilde\theta_0-\theta_0\)\) \right\rvert,
\end{equation*}
for $K$ sufficiently large. We now have
\begin{align*}
  \left\lvert \theta_0 - 
      \widetilde\theta_0\right\rvert &\gtrsim \left\lvert
      |\alpha_0|^{2\si}-|\widetilde 
    \alpha_0|^{2\si} \right\rvert
+|\alpha_1|^{2\si-2}\left\lvert  |\alpha_0|^{2}-|\widetilde
    \alpha_0|^{2} \right\rvert\\
&\gtrsim \delta + \(\rho K^{|s|}\)^{2\si-2}\delta. 
\end{align*}
Now we see that if we assume $\si\ge 2$, the left hand side can be
estimated from below by $1/\delta$, provided $N$ is sufficiently large,
and we conclude like for the first point. 

To prove the last point in Theorem~\ref{theo:instab}, we resume the
argument of \cite{Molinet}. Fix $\alpha_0\in \C\setminus\{0\}$, and let
$\alpha_1 \in \C$ to be fixed later. As $K\to
\infty$, we have:
\begin{equation*}
  u(0,\cdot)\rightharpoonup \alpha_0=:\underline u(0,\cdot)\text{
    weakly in }L^2(\T^d)\quad 
  ;\quad \|u(0)\|_{L^2}^2\to  |\alpha_0|^2 +|\alpha_1|^2. 
\end{equation*}
For any $t>0$, we have, as $K\to \infty$,
\begin{equation*}
  u(t,x)\rightharpoonup \alpha_0 e^{-i\l t \theta_0} \text{ weakly in
  }L^2(\T^d), 
\end{equation*}
where
\begin{equation*}
  \theta_0 = \sum_{n=0}^\si \(
    \begin{array}[c]{c}
      \si+1 \\
n
    \end{array}\)
\(
    \begin{array}[c]{c}
      \si \\
n
    \end{array}\)|\alpha_0|^{2\si-2n}|\alpha_1|^{2n}.
\end{equation*}
Note that for any $\alpha_0\in \C\setminus\{0\}$ and any angle
$\theta\in [0,2\pi[$, we can find $\alpha_1\in \C$ so that
$\theta_0=\theta+ |\alpha_0|^{2\si}$. On the other hand, the solution
to \eqref{eq:nlsT} with initial data $\alpha_0$ is given by 
\begin{equation*}
  \underline u(t,x)=\alpha_0 e^{-i\l t|\alpha_0|^{2\si}}. 
\end{equation*}
We infer
\begin{equation*}
  w-\lim_{N\to \infty} u(t,x) - \underline u(t,x) = \alpha_0e^{-i\l
    t|\alpha_0|^{2\si}} \( e^{-i\l t \theta} -1\). 
\end{equation*}
For all $t\not =0$, one can then choose $\theta$ so that $\l t\theta
\not\in2\pi\Z$. The discontinuity at $\alpha_0$ of the map
$\alpha_0\mapsto \underline u(t)$, from $L^2(\T^d)$ equipped with its
weak topology into $\(C^\infty(\T^d)\)^*$, follows. 
\end{proof}

\appendix

\section{A more general set of initial phases}
\label{sec:appendix}

We can actually replace Assumption~\ref{as:nosmalldiv} with the following more
general one:
\begin{assumption}
\label{as:nosmalldivgen}
There exist $b\ge0$, $c>0$ such that for all $(\ell_1,\dots,\ell_{2\si+1}) 
\in N $, 
\begin{equation*}
  \delta\(\ell_1,\dots,\ell_{2\si+1}\) \equiv
\left\lvert \left\lvert
    \sum_{p=1}^{2\si+1}(-1)^{p+1}\k_{\ell_p}\right\rvert^2
  -\sum_{p=1}^{2\si+1}(-1)^{p+1}\left\lvert \k_{\ell_p}\right\rvert^2
\right\rvert  
\end{equation*}
satisfies:
$$\delta(\ell_1,\dots,\ell_{2\si+1}) \ge c \<\k_{\ell_1}\>^{-b} 
\dots \<\k_{\ell_{2\si+1}}\>^{-b} .$$
\end{assumption}
In \S\ref{sec:justif}, we have considered the case $b=0$. 
However, allowing constants $b>0$, we show that the assumption is 
satisfied by wave vector sets included in generic finitely generated 
nets.
\begin{proposition}
For all $p\in\N^*$, there exist $C,b > 0$ and $Z\subset\R^{dp}$ with zero 
Lebesgue measure such that, for all
$(\k_1,\dots,\k_p)\in\R^{dp}\setminus Z$,  
the set $(\k_j)_{j\in J}$ constructed from these initial wave vectors 
$\{\k_j\}_{j\in J_0}$ satisfies Assumption~\ref{as:nosmalldivgen}.
\end{proposition}
\begin{proof}
We shall prove that the above result holds when
Assumption~\ref{as:nosmalldivgen} is replaced by the stronger one,
where $N$ is replaced by $J^{2\si+1}$. 

All the wave vectors we consider belong to the group generated by 
$\{\k_1,\dots,\k_p\}$. Thus, to each $\ell_k\in J$ corresponds 
$(\alpha_{k,1},\dots,\alpha_{k,p})\in\Z^p$, such that:
$\k_{\ell_k} = \alpha_{k,1}\k_1 + \dots + \alpha_{k,p}\k_p$. 
With this notation, for all $(\ell_1,\dots,\ell_{2\si+1}) \in J^{2\si+1}$, 
we have:
\begin{align*}
\delta(\ell_1,\dots,\ell_{2\si+1}) 
& \quad = \left| 
\Big| 
\sum_{k=1}^{2\si+1} (-1)^{k+1} \sum_{j=1}^p \alpha_{k,j} \k_j 
\Big|^2 
+ \sum_{m=1}^{2\si+1} (-1)^{m+1} \Big| \sum_{j=1}^p \alpha_{m,j} 
\k_j \Big|^2
\right| \\
& \quad = \left| 
\sum_{i,j=1}^p \left( \sum_{k,\ell=1}^{2\si+1} (-1)^{k+\ell} \alpha_{k,i}
\alpha_{\ell,j} 
- \sum_{m=1}^{2\si+1} (-1)^m \alpha_{m,i} \alpha_{m,j} \right)
\k_i \cdot \k_j 
\right| .
\end{align*}
Now, a standard Diophantine result (see \emph{e.g.}
\cite{AlinhacGerardUS,Craig00}) 
ensures that,  
for all choice of $(\k_i\cdot\k_j)_{1\le i,j\le p}$ but in some subset 
of $\R^{p^2}$ with measure zero, we have, for some $b'\ge0$ and $C'>0$:
\begin{equation*} \label{diophbound}
\forall (\beta_{i,j})_{1\le i,j\le p}\in\Z^{p^2}\setminus\{0\}, \quad 
\left| \sum_{i,j=1}^p \beta_{i,j} \k_i\cdot\k_j \right| \ge C' 
\left( \sum_{i,j=1}^p \left| \beta_{i,j} \right| \right)^{-b'}. 
\end{equation*}
Such an estimate is then valid for almost all $(\k_1,\dots,\k_p)$ in 
${(\R^d)}^p$. We apply it with 
$$\beta_{i,j} = \sum_{k,\ell=1}^{2\si+1} (-1)^{k+\ell} \alpha_{k,i}
\alpha_{\ell,j} - \sum_{m=1}^{2\si+1} (-1)^m \alpha_{m,i} \alpha_{m,j},$$
so that
\begin{align*}
\sum_{i,j=1}^p |\beta_{i,j}|
& \le 2 \sum_{k,\ell=1}^{2\si+1} |\alpha_{k,\cdot}| |\alpha_{\ell,\cdot}| \\
& \le 2(2\si+1)^2 \prod_{k=1}^p \<\alpha_{k,\cdot}\>^2. 
\end{align*}
Now, choosing $\k_1,\dots,\k_p$ $\Q$-linearly independent (which is true 
almost surely), we get that there exists a constant $c>0$ such that 
$$\forall \alpha \in \Q^p, \quad 
|\alpha_1| + \dots + |\alpha_p| \le c \sum_{j=1}^d |(\alpha_1 \k_1 + 
\dots + \alpha_p \k_p)_j|.$$
Increasing $c$ if necessary, so that $c\ge1$, we get, when 
$\k_{\ell_k} = \alpha_{k,1}\k_1 + \dots + \alpha_{k,p}\k_p$: 
$\<\alpha_{k,\cdot}\> \le c \<\k_{\ell_k}\>$. 
Finally, using the constants $b'$ and $C'$ from above,
the desired estimate follows with $b=2b'$ and  $C=(2(2\si+1)^2c^2)^{-b'}C'$.
\end{proof}

Under Assumption~\ref{as:nosmalldivgen} (which is fairly  general for plane
waves, in view of the above proposition), we can easily adapt the
analysis of \S\ref{sec:justif}. Essentially, we have to (possibly)
strengthen the assumptions on the initial profile, in the case of
$\mathcal M = \R^d$, where we generalize Assumption \ref{as:regularity} to:
\begin{assumption}
On $\mathcal M = \R^d$, the initial amplitudes satisfy:
\label{as:regularitygen}
\begin{equation*}  
  \begin{aligned}
   \forall |\beta|\le 2, & \
(\<\k_j\>^b\partial^\beta_x\alpha_j)_{j\in J_0}\in E(\R^d),\\ 
 \forall |\beta|\le 1, &\ (\<\k_j\>^{1+b}\partial^\beta_x\alpha_j)_{j\in J_0}
\in E(\R^d).
  \end{aligned}
\end{equation*}
\end{assumption}

From Proposition \ref{prop:existprof}, these data produce a solution
$(a_j)_{j\in J}\in C([0,T],E(\mathcal M))$ to the profile system
\eqref{eq:transportsystemint}. We consequently define the approximate solution $u_{\rm app}^\eps$ as before
\begin{equation*}
  u_{\rm app}^\eps(t,x) = \sum_{j\in J}a^\eps_j(t,x) e^{i\phi_j(t,x)/\eps},
\end{equation*}
where the sequence $(a_j)_{j\in J}$ is now such that
\begin{align*}
  &\(\<\k_j\>^b\partial^\beta_xa_j\)_{j\in J}\in C([0,T],E(\mathcal M)),
  \quad |\beta|\le 2,\\
& \(\<\k_j\>^{1+b} \partial^\beta_xa_j\)_{j\in J}\in
C([0,T],E(\mathcal M)), \quad |\beta|\le 1.
\end{align*}
We can then reproduce the analysis of \S\ref{sec:justif}:
Lemma~\ref{lem:Rest} is still valid under 
Assumption~\ref{as:nosmalldivgen} and \ref{as:regularitygen}, by
straightforward verification. Then one just has to notice that this is
the only step where the absence of small divisors plays a role in the
proof of Theorem~\ref{theo:justif}. Therefore,
Theorem~\ref{theo:justif} remains valid under
Assumption~\ref{as:nosmalldivgen} and \ref{as:regularitygen}.

\bibliographystyle{amsplain}
\bibliography{wnlgo}

\end{document}